\documentclass[11pt,reqno]{amsart}
\usepackage[%
   pdftitle={{On general transport equations with
abstract boundary conditions. The case of divergence free force
field.}},
   pdfauthor={Luisa Arlotti, Jacek Banasiak and Bertrand Lods},
   pdfsubject={{Semigroups for general transport equations with abstract
boundary conditions}},
   pdfkeywords={Transport equation, Boundary
conditions, $C_0$-semigroups, Characteristic curves},
  colorlinks=true,
  linkcolor=red,
  citecolor=blue,
]{hyperref}
\usepackage[T1]{fontenc}
\usepackage{mathrsfs}
\usepackage{amsmath,amsthm,amssymb}
\usepackage{amscd,indentfirst,epsfig}

\usepackage{latexsym}
\usepackage{times}
\usepackage{enumerate}
\usepackage{mathrsfs}
\usepackage{stmaryrd}
\usepackage{amsopn}
\usepackage{amsmath}
\usepackage{amssymb}
\usepackage{amsfonts}
\usepackage{amsbsy}
\usepackage{amscd,indentfirst}
\usepackage{hyperref}
\usepackage{amsfonts,amsmath,latexsym,amssymb,verbatim,amsbsy}
\usepackage{amsthm}
\usepackage{colordvi}
\renewcommand{\theequation}{\arabic{section}.\arabic{equation}}

\newtheorem{theo}{{Theorem}} [section]
\newtheorem{lemme}[theo]{{Lemma}}
\newtheorem{propo}[theo]{{Proposition}}
\newtheorem{cor}[theo]{{Corollary}}
\newtheorem{hyp}{{Assumption}}
\newtheorem{nb}[theo]{{Remark}}
\newtheorem{defi}[theo]{{Definition}}
\theoremstyle{definition}
\newtheorem{exa}[theo]{{Example}}

\newenvironment{proof1}{\noindent {\textsc{Proof of Proposition \ref{moll2}.}}}{\hfill $\square$\bigskip}
\newenvironment{preuve2}{\noindent {\textsc{Proof of Theorem \ref{representation}.}}}{\hfill $\square$\bigskip}
\def \leq {\leqslant}
\def \geq {\geqslant}
\numberwithin{equation}{section}
\def\ind#1{\lower5pt\hbox{$\scriptstyle #1$}}

\def \O {\mathbf{\Omega}}
\newcommand{\Con}{\ensuremath{\mathscr{C}}}

\newcommand{\Drond}{\ensuremath{\mathscr{D}}}
\def \x {\mathbf{x}}
\def \z {\mathbf{z}}
\def\y {\mathbf{y}}
\def\lp {L^1_+}
\def\lm {L^1_-}
\def \ff {\mathscr{F}}

\def \d {\mathrm{d}}
\def \D {\Drond}

\def \e {\epsilon}
\def \F {\mathbf{F}}
\def \A {\mathcal{A}}
\def \uot {(U_0(t))_{t \geq 0}}

\def \ds {\displaystyle}
\def \T {\mathcal{T}}
\def \B {\mathsf{B}}

\def \t {\tau}
\title[A new approach to transport equations]{A new approach to transport equations associated to a regular field: trace results and
well-posedness.}

\author{L. \textsc{Arlotti}, J. \textsc{Banasiak} \& B. \textsc{ Lods}}

\address{\newline \noindent \textit{\textbf{Luisa Arlotti}} \newline \noindent Dipartimento di Ingegneria Civile, Universit\`a di
Udine, via delle Scienze 208,\newline 33100 Udine, Italy.\newline
{\tt luisa.arlotti@uniud.it}}

\address{\newline \noindent \textit{\textbf{Jacek Banasiak}} \newline School of
Mathematical Sciences,
 University of KwaZulu--Natal,\newline Durban 4041, South Africa. \newline
 {\tt banasiak@ukzn.ac.za}}

\address{\newline \noindent \textit{\textbf{Bertrand Lods}}
\newline
\noindent Laboratoire de Math\'{e}matiques, CNRS UMR 6620, Universit\'{e}
Blaise Pascal (Clermont-Ferrand 2), 63177 Aubi\`{e}re Cedex, France.
\newline
\noindent{\tt bertrand.lods@math.univ-bpclermont.fr} }
\begin{document}
 \maketitle
\begin{abstract}
We generalize known results on transport equations associated to a
Lipschitz field $\ff$ on some subspace of $\mathbb{R}^N$ endowed
with some general space measure $\mu$. We provide a new definition
of both the transport operator and the trace measures over the
incoming and outgoing parts of $\partial \O$ generalizing known
results from \cite{beals, voigt}. We also prove the well-posedness
of some suitable boundary-value transport problems and describe in
full generality the generator of the transport semigroup with
no-incoming boundary conditions.
\end{abstract}
\medskip
%

%
\section{Introduction}
\noindent

In this paper we present new methodological tools to investigate the
  well-posedness of the general transport equation
\begin{subequations}\label{1}
\begin{equation}\label{1a}
\partial_t f(\x,t)+\ff(\x)\cdot \nabla_\x
f(\x,t) =0 \qquad (\x \in\O, \:t
> 0),\end{equation} supplemented by boundary condition
\begin{equation}\label{1b}
f_{|\Gamma_-}(\y,t)=\psi_-(\y,t), \qquad \qquad (\y \in \Gamma_-, t
>0),
\end{equation}
and the initial condition
\begin{equation}\label{1c}f(\x,0)=f_0(\x), \qquad \qquad (\x \in \O). \end{equation}\end{subequations}
Here $\O$ is a sufficiently smooth open subset of $\mathbb{R}^N$,
$\Gamma_{\pm}$ are suitable boundaries of the phase space and $\psi_-$
is a given function of the trace space $L^1(\Gamma_-,\d\mu_-)$
corresponding to the boundary $\Gamma_-$ (see Section 2 for
details).

The present paper is part of a series of papers on transport
equations with general vector fields \cite{m2as,abl} and introduce
all the methodological tools that  allow us not only to solve the
initial-boundary problem \eqref{1} but also to treat in \cite{abl}
  the case of abstract boundary
conditions relying the incoming and outgoing fluxes, generalizing
the results of \cite{beals}.

The main novelty of our approach is that we assume $\mathbb{R}^N$ to
be endowed with a general positive Radon measure $\mu$. Here by a
Radon measure we understand a Borel measure (or its completions, see
\cite[p.~332]{Roy}) which is finite on compact sets. As we shall see
it further on, taking into account such general Radon measure $\mu$
leads to a large amount of technical difficulties, in particular in
the definition of trace spaces and in the derivation of Green's
formula. Moreover, for such a measure $\mu$, it is far from being
trivial to identify the vector field $\ff \cdot \nabla_x$ to the
time derivative along the characteristic curves (as done in
\cite[Formulae (5.4) \& (5.5), p.392]{beals}): the main difficulty
stemming from the impossibility of applying classical convolution
arguments (and the so-called Friedrich's lemma). We overcome this
difficulty by introducing new mollification techniques  along the
characteristic curves. Let us explain in more details our general
assumptions:

\subsection{General assumption and motivations}
The transport coefficient $\ff$ is a \textit{time independent}
vector field
$\ff \::\:\mathbb{R}^N  \longrightarrow
\mathbb{R}^N$ which is (globally) \textit{Lipschitz-continuous} with Lipschitz
constant $\kappa > 0$, i.e.
\begin{equation}\label{lipsc}
|\ff(\x_1)-\ff(\x_2)|\leq \kappa |\x_1-\x_2| \qquad \text{ for any
} \:\x_1,\x_2 \in \mathbb{R}^N.\end{equation}
 Clearly, one can associate a flow $(T_t)_{t \in \mathbb{R}} $
to this field $\ff$ (with the notations of Section \ref{sub:chara},
 $T_t=\mathbf{\Theta}(\cdot, t, 0)$)
and we make the following fundamental assumption (known as
\textit{\textbf{Liouville's Theorem}} whenever $\mu$ is the Lebesgue
measure) on $\ff$:
\begin{hyp}\label{ass:h2} The measure $\mu$ is invariant under the
flow $(T_t)_{t \in \mathbb{R}}$, i.e. $\mu(T_t A) = \mu(A)$ for any
measurable subset  $A \subset \mathbb{R}^N$ and any $t \in
\mathbb{R}$.\end{hyp}

\begin{nb}  Notice that, whenever $\mu$ is the Lebesgue measure over
$\mathbb{R}^N$, it is well-known that Assumption \ref{ass:h2} is
equivalent to $ \mathrm{div}(\ff(\x))=0$  for any $\x \in
\mathbb{R}^N.$ More generally, by virtue of \cite[Remark 3 \&
Proposition 4]{ambrosio}, Assumption \ref{ass:h2} holds for a
general
 Borel measure $\mu$ provided the  field $\ff$ is locally
 integrable with respect to $\mu$ and
\textbf{\textit{divergence-free}} with respect to $\mu$ in the sense
that
$$\int_{\mathbb{R}^N}\ff(T_t(\x)) \cdot
\nabla_{\x}f(T_t(\x))\d\mu(\x)=0, \qquad \qquad \forall t \in
\mathbb{R}$$ for any infinitely differentiable function $f$ with
compact support.
 \end{nb}

\noindent  A typical example of such a transport equation is the
so-called Vlasov equation for which:
\begin{enumerate}[i)\:] \item The phase space $\O$ is given by the
cylindrical domain $\O=\mathcal{D}\times \mathbb{R}^3 \subset
\mathbb{R}^6$ where $\mathcal{D }$ is a sufficiently
 smooth open subset of $\mathbb{R}^3$, referred to as the
\textit{position space}, while the so--called \textit{velocity
space} is here given by $\mathbb{R}^3.$ The measure $\mu$ is given
by $\d\mu(\x)=\d x \d\beta(v)$ where $\beta $ is a suitable Radon
measure on $\mathbb{R}^3$, e.g. Lebesgue measure over $\mathbb{R}^3$
for continuous models or combination of Lebesgue measures over
suitable spheres for the multigroup model. \item For any $\x=(x,v)
\in \mathcal{D} \times \mathbb{R}^3$,
\begin{equation}\label{vlasov}\ff(\x)=(v,\F(x,v)) \in \mathbb{R}^6\end{equation} where
$\F=(\F_1,\F_2,\F_3)$ is a time independent force field over
$\mathcal{D} \times \mathbb{R}^3$ satisfying  Assumption
\ref{ass:h2} and (\ref{lipsc}). The free transport case,
investigated in \cite{voigt,labl}, corresponds to $\F=0$.
\end{enumerate}


The existence of solution to the transport equation \eqref{1a} is a
classical matter when considering the whole space $\O=\mathbb{R}^N$.
In particular, the concept of renormalized solutions allows to
consider irregular transport coefficient $\ff(\cdot)$ (see
\cite{diperna} and the recent contributions \cite{ambrosio,lions})
which is of particular relevance in fluid mechanics.

On the other hand, there are few results addressing the
initial-boundary value problem \eqref{1}, possibly due to
difficulties created by the boundary conditions \eqref{1b}. We
mention here the seminal works by C. Bardos \cite{bardos}, and by R.
Beals and V. Protopopescu \cite{beals} (see also \cite{proto,vdm}).
Let us however mention that the results of \cite{beals,proto}
introduce restrictive assumptions on the characteristics of the
equation. For instance, fields with 'too many' periodic trajectories
create serious difficulties. They are however covered in a natural
way by the theory presented here, see Examples \ref{exa1} \&
\ref{exa2}.

\subsection{Presentation of the results}

In this paper, we revisit and generalize the afore-mentioned results
to the general case $\F\neq 0$ and for a general Radon measure
$\mu$.  The latter, in particular, leads to numerous technical
problems such as e.g. determination of suitable measures $\mu_{\pm}$
over the `incoming' and `outgoing' parts $\Gamma_{\pm}$ of $\partial
\O$. We provide here a general construction of these `trace
measures' generalizing, and making more precise, the results of
\cite{beals,proto}. This construction allows us to establish
Proposition \ref{prointegra} which allows to compute integrals over
$\mathbf{\Omega}$ via integration along the integral curves of
$\ff(\cdot)$ coming from the boundary $\partial \mathbf{\Omega}$,
and which is free from some restrictive assumptions of \textit{op.
cit}. In particular, we present a new proof of the Green formula
clarifying and removing gaps of the proofs in \cite{beals,proto}. Of
course, the boundary condition \eqref{1b} we treat here is less
general than the abstract ones investigated in \cite{beals,proto}
but, as we already mentioned  it, the tools we introduce here will
allow us to generalize, in a subsequent paper \cite{abl}, the
results of the \textit{op. cited} by dealing with abstract boundary
conditions.

 Another major difficulty, when dealing with a
general Radon measure $\mu$, is to provide a precise definition of
the transport operator $\T_\mathrm{max}$ associated to \eqref{1}. It
appears quite natural to define the transport operator
$\T_\mathrm{max}$ (with its maximal domain on $L^1(\O,\d\mu)$) as a
\textit{{weak directional derivative along the characteristic
curves}} in the $L^1$-sense. However, it is not clear {\it a priori}
that any function $f$ for which  the weak directional derivative
exists in $L^1(\O,\d\mu)$ (with appropriate and minimal class of
test-functions)  admits a trace over $\Gamma_{\pm}$. With the aim of
proving  such a trace result, we provide here a new characterization
of the transport operator related to a \textit{{mild
representation}} of the solution to \eqref{1}. Namely, we prove
(Theorem \ref{representation}) that the domain $\D(\T_\mathrm{max})$
(as defined in Section 3), is precisely the set of functions $f \in
L^1(\O,\d\mu)$ that admits a representative which is
\textit{absolutely continuous along almost any characteristic
curve}.

Note that in the classical case when $\mu$ is the Lebesgue measure,
such a representation is known to be true \cite[Appendix]{diperna}.
Actually, in this case, one defines the domain $\D(\T_\mathrm{max})$
 as the set of all $f \in L^1(\O,\d\mu)$ for which
the directional derivative $-\ff \cdot \nabla f$ exists in the
distributional sense and belongs to $L^1(\O,\d\mu)$. Then, by
convolution arguments, it is well-known that the set $\Con_0^1(\O)
\cap \D(\T_\mathrm{max})$ is dense in $\D(\T_\mathrm{max})$ for the
graph norm $\|f\|=\|f\|+\|\ff \cdot \nabla f\|.$

The question is much more delicate for a   general Radon measure
$\mu$. Indeed, in such a case, the convolution argument used in the
case of the Lebesgue measure does not apply anymore. Our strategy to
prove the characterization of $\T_\mathrm{max}$ is also based on a
convolution argument but it uses  \textit{mollification technique
along the characteristic curves} as developed in Section 3. Such a
result shall allow us to obtain a rigorous derivation of Green's
formula, clarifying some results of \cite{beals}.
\subsection{Plan of the paper} The organization of the paper is as follows. In Section 2
we introduce main tools used throughout the paper and present the
aforementioned new results concerning \textit{integration over the
characteristic curves} of $\ff$ as well as \textit{a new
construction of the boundary measures} over the `incoming' and
`outgoing' parts $\Gamma_{\pm}$ of $\partial \O$ which generalizes and
clarifies that of \cite{beals,proto}. In Section 3 we provide a
construction of the maximal transport operator $\T_\mathrm{max}$. It
is defined in a weak sense, through its action on suitably defined
test functions. The fundamental result of this section shows that
any function in the domain $\D(\T_\mathrm{max})$ admits a
representation which is absolutely continuous along almost any
characteristic which, in turn, allows for existence of its traces on
the outgoing and incoming parts of the boundary. In Section 4 we
apply the results of Section 3 to prove well-posedness of the
time--dependent transport problem with no reentry boundary
conditions associated with $\T_\mathrm{max}$. Moreover, we consider
the corresponding stationary problem and, as a by-product, we
recover  a new proof of the \textit{Green formula}.

\section{Integration along the characteristics}
\subsection{\textbf{Characteristic curves}}\label{sub:chara}   A crucial role in our study
is played  by the characteristic curves associated to the
field $\ff$.
Precisely, for any $\x \in \mathbb{R}^N$ and $t \in \mathbb{R}$,
consider the initial-value problem
\begin{equation}\label{chara}
\begin{cases}
\dfrac{\d }{\d s}\mathbf{X}(s)=\ff(\mathbf{X}(s)), \qquad (s \in \mathbb{R});\\
\mathbf{X}(t)=\x.
\end{cases}\end{equation}
Since $\ff$ is Lipschitz continuous on $\mathbb{R}^N$, Eq.
\eqref{chara} has a unique {\it global in time} solution and this
allows to define the flow--mapping $\mathbf{\Theta}$ : $\mathbb{R}^N \times
\mathbb{R}\times \mathbb{R} \to \mathbb{R}^N$, such that, for
$(\x,t) \in \mathbb{R}^N \times \mathbb{R}$, the mapping:
$$\mathbf{X}(\cdot)\::\:s \in \mathbb
{R} \longmapsto \mathbf{\Theta}(\x,t,s)$$ is the only solution of
Eq. \eqref{chara}.  
Being concerned with solutions to the transport equation \eqref{1}
in the region $\O$, we have to introduce the definition of stay
times of the characteristic curves in $\O$:
\begin{defi} For any $\x \in \O$, define
$\tau_{\pm}(\x)=\inf \{s > 0\,;{\mathbf{\Theta}}(\x,0,\pm s) \notin \O\},$
with the convention that $\inf \varnothing=\infty,$ and set $\t(\x)
= \t_+(\x) + \t_-(\x). $\end{defi} In other words, given $\x \in
\O$, $I_{\x}=(-\t_-(\x),\t_+(\x))$ is the maximal interval for which
${\mathbf{\Theta}}(\x,0,s)$ lies in $\O$ for any $s \in I_{\x}$ and
$\t(\x)$ is the length of the interval $I_{\x}$. Notice that $0\leq
\t_{\pm}(\x) \leq \infty$. Thus, the function ${\mathbf{\Theta}}$
restricted to the set
$$\mathbf{\Lambda}:=\bigg\{(\x,t,s)\,;\,\x \in \O,\,t \in
\mathbb{R}\,,\,s \in \left(t-\t_-(\x),t+\t_+(\x)\right)\;\bigg\}$$
is such that ${\mathbf{\Theta}}(\mathbf{\Lambda})=\O$.  Note that
here we \textit{\textbf{do not}} assume that the length of the
interval $I_{\x}=(-\t_-(\x),\,\t_+(\x))$ is
\textit{\textbf{finite}}. In particular, $I_{\x}=\mathbb{R}$ for any
stationary point $\x$ of $\ff$, i.e. $\ff(\x)=0$. If $\tau(\x)$ is
finite, then the function $\mathbf{X}\::\:s \in I_{\x} \longmapsto
{\mathbf{\Theta}}(\x,0,s)$ is bounded since $\ff$ is Lipschitz
continuous. Moreover, still by virtue of the Lipschitz continuity of
$\ff$, the only case when $\t_{\pm}(\x)$ is finite is when
${\mathbf{\Theta}}(\x,0,\pm s)$ reaches the boundary $\partial\O$ so
that ${\mathbf{\Theta}}(\x,0,\pm \tau_{\pm}(\x)) \in
\partial\O$. We note that, since $\ff$ is Lipschitz
around each point of $\partial\O$, the points of the set $\{\y \in
\partial\O\,;\,\ff(\y)=0\}$ (introduced in
\cite{beals,proto}) are equilibrium points of the $\ff$ and
cannot be reached in finite time.
 \begin{nb}
We emphasize that periodic trajectories which do not meet the
boundaries have $\tau_\pm =\infty$ and thus are treated as infinite
though geometrically they are bounded.
\end{nb}

Finally we mention that it is not difficult to prove that the
mappings $\tau_{\pm}$ : $\O \to \mathbb{R}^+$ are lower
semicontinuous and therefore measurable, see e.g., \cite[p.\
301]{arloban}

The flow ${\mathbf{\Theta}}(\x,t,s)$ defines, at each instant $t$, a
mapping of the phase space $\O$ into $\mathbb{R}^N$. Through this
mapping, to each point $\x$ there corresponds the point
$\x_{s,t}={\mathbf{\Theta}}(\x,t,s)$ reached at time $s$ by
the point which was at $\x$ at the `initial' time $t$. 
The flow ${\mathbf{\Theta}}$, restricted to $\mathbf{\Lambda}$, has
the properties:
\begin{propo}\label{Phiprop} Let $\x \in \O$ and $t \in \mathbb{R}$ be fixed. Then,
\begin{enumerate}[(i)\:]
\item ${\mathbf{\Theta}}(\x,t,t)=\x.$ \item
${\mathbf{\Theta}}({\mathbf{\Theta}}(\x,t,s_1),s_1,s_2)={\mathbf{\Theta}}(\x,t,s_2),
\quad \forall  s_1,s_2 \in (t-\t_-(\x),t+\t_+(\x)).$ \item
${\mathbf{\Theta}}(\x,t,s)={\mathbf{\Theta}}(\x,t-s,0)={\mathbf{\Theta}}(\x,0,s-t),
\qquad \forall s \in (t-\t_-(\x),t+\t_+(\x)).$ \item
$\left|{\mathbf{\Theta}}(\x_1,t,s)-{\mathbf{\Theta}}(\x_2,t,s)\right|
\leq \exp(\kappa |t-s|)|\x_1-\x_2|$ for any $\x_1,\x_2 \in \O$, $s -
t \in I_{\x_1} \cap I_{\x_2}.$
\end{enumerate}\end{propo}


An important consequence of $(iii)$ above is that
${{\mathbf{\Theta}}}(\x,0,s)={\mathbf{\Theta}}(\x,-s,0)$ for any $\x
\in\O$, $0 \leq s \leq \tau_+(\x).$ Therefore, from now on, to
shorten notations we shall denote
$${\mathbf{\Phi}}(\x,t)={\mathbf{\Theta}}(\x,0,t), \qquad \forall t \in \mathbb{R},$$
so that ${\mathbf{\Phi}}(\x,-t)={\mathbf{\Theta}}(\x,t,0)$, $t \in
\mathbb{R}.$ We define the incoming and outgoing part of the
boundary $\partial \O$ through the flow ${\mathbf{\Phi}}$:
\begin{defi} The
\textit{incoming} $\Gamma_{-}$ and the \textit{outgoing} $\Gamma_+$
parts of the boundary $\partial \O$ are defined by:
\begin{equation}\label{gammapm}
\Gamma_{\pm}:=\left\{\y \in \partial \O\,;\exists \x \in \O,\, \t_{\pm}(\x)
< \infty \text{ and } \y={\mathbf{\Phi}}(\x,\pm
\t_{\pm}(\x))\,\right\}.\end{equation}
\end{defi}
Properties of ${\mathbf{\Phi}}$ and of $\t_\pm$ imply that
$\Gamma_{\pm}$ are Borel sets. It is possible to extend the definition
of $\t_{\pm}$ to $\Gamma_{\pm}$ as follows. If $\x\in\Gamma_-$ then we put
$\t_-(\x) = 0$ and denote $\t_+(\x)$ the length of the integral
curve having $\x$ as its left end--point; similarly if
$\x\in\Gamma_+$ then we put $\t_+(\x) = 0$ and denote $\t_-(\x)$ the
length of the integral curve having $\x$ as its right endpoint. Note
that this definition implies that $\tau_{\pm}$ are measurable over $\O
\cup \Gamma_- \cup \Gamma_+$.



Let us illustrate the above definition of $\Gamma_{\pm}$ by two
simple $2D$ examples:
\begin{exa}[\textit{\textbf{Harmonic oscillator in a rectangle}}]\label{exa1} Let $\O=(-a,a) \times (-\xi,\xi)$ with $a,\xi >0$ and
let us consider the harmonic oscillator force field
\begin{equation}\label{harmonic}
\ff(\x)=(v,-\omega^2 x), \qquad \text{ for any } \x=(x,v) \in
\O\end{equation} where $\omega
>0$. We take as $\mu$ the Lebesgue measure over $\mathbb{R}^2$
and, since $\ff$ is divergence-free, Assumption
\ref{ass:h2} is fulfilled. In this case, for any $\x_0=(x_0,v_0)
\in \O$, the solution $(x(t),v(t))={\mathbf{\Phi}}(\x_0,t)$ to the
characteristic equation $\frac{\d}{\d
t}\mathbf{X}(t)=\ff\left(\mathbf{ X}(t)\right),$ $
\mathbf{X}(0)=\x_0$, given by
$${\mathbf{\Phi}}(\x_0,t)=\left(x_0\cos(\omega t)+\frac{v_0}{\omega}\sin(\omega
t)\,;\,-x_0\omega \sin(\omega t)+v_0 \cos(\omega t)\right),$$ is
such that
$$\omega^2 x^2(t)+v^2(t)=\omega^2 x_0^2+v_0^2, \qquad t \in (-\t_-(\x_0),\t_+(\x_0))$$ which means that the
integral curves associated to $\ff$ are \textit{ellipses} centered
at $(0,0)$ and oriented in the counterclockwise direction. Now,
$$\partial \O=\bigg(\{-a\} \times [-\xi,\xi]\bigg) \bigcup \bigg(\{a\} \times
[-\xi,\xi]\bigg) \bigcup \bigg([-a,a] \times \{-\xi\}\bigg) \bigcup
\bigg([-a,a] \times \{ \xi\}\bigg)$$ and it is easy to check that
$$\Gamma_{\pm}=\bigg(\{\pm a\} \times (-\xi,0]\bigg)
\bigcup \bigg(\{\mp a\} \times [0,\xi)\bigg) \bigcup \bigg([0,a)
\times \{\pm \xi\}\bigg) \bigcup \bigg((-a,0] \times \{\mp
\xi\}\bigg).$$ Notice that $\Gamma_{+} \cap
\Gamma_-=\{(a,0),\,(0,\xi),\,(-a,0),\,(0,-\xi)\}$ and $$\partial \O
\setminus \left(\Gamma_+ \cup
\Gamma_-\right)=\{(a,\xi),(a,-\xi),(-a,\xi),(-a,\xi)\}$$ is a
discrete  set (of linear Lebesgue measure zero).
\end{exa}
\begin{exa}[\textit{\textbf{Hamonic oscillator in a
stadium}}]\label{exa2} Consider now the two-dimensional phase space
(where $\mathbb{R}^2$ is still endowed with the Lebesgue measure
$\mu$):
$${\O}=\{\x=(x,v) \in \mathbb{R}^2\,;\,x^2+v^2 < 2 \:\text{
and } -1 < v < 1\}$$ and consider the harmonic oscillator force
field $\ff$ given by \eqref{harmonic} with $\omega=1$ for
simplicity. Then, the integral curves associated to $\ff$ are
 \textit{circles} centered at $(0,0)$ and oriented in the counterclockwise
direction. In this case, one can see that
$$\Gamma_{\pm}=\{(x,-1)\,;\,-1 < \pm\, x \leq 0\} \cup \{(x,1)\,;\,0\leq \pm
\,x<1\}.
$$
In particular, one sees that $\partial \O \setminus \big(\Gamma_+
\cup \Gamma_-\big)=\left\{(x,v)\in \mathbb{R}^2\,;\,x^2+v^2=2\,;\,
-1 \leq v \leq 1\right\}$ is a 'big' part of the boundary $\partial
\O$ (with positive linear Lebesgue measure). Notice also that
$\t_+(\x)=+\infty$ for any $\x=(x,v)$ with $x^2+v^2 <1$.\end{exa}
The main aim of the present discussion is to represent $\O$ as a
collection of characteristics running between points of $\Gamma_-$
and $\Gamma_+$ so that the integral over $\O$ can be split into
integrals over $\Gamma_-$ (or $\Gamma_+$) and along the
characteristics. However, at present we cannot do this in a precise
way since, in general, the sets $\Gamma_+$ and $\Gamma_- $ do not
provide a partition of $\partial \O$ as there may be `too many'
characteristics which extend to infinity on either side. Since we
have not assumed $\O$ to be bounded, $\Gamma_-$ or $\Gamma_+$ may be
empty and also we may have characteristics running from $-\infty$ to
$+\infty$ such as periodic ones. Thus, in general, characteristics
starting from $\Gamma_-$ or ending at $\Gamma_+$ would not fill the
whole $\O$ and, to proceed, we have to construct an auxiliary set by
extending $\O$ into the time domain and use the approach of
\cite{beals} which is explained below.

\subsection{\textbf{Integration along characteristics}} For any $0 < T < \infty$, we define the domain
$$\O_T=\O \times (0,T)$$
and the measure $\d\mu_T=\d\mu\otimes\d t$ on $\O_T$.
Consider the vector field over $\O_T$:
$$Y={\partial_t} + \ff(\x)\cdot \nabla_\x=\mathscr{A}(\xi)\cdot \nabla_{\xi}$$
where $\mathscr{A}(\xi)=(\ff(\x),1)$ for any $\xi=(\x,t)$. We can
define the characteristic curves of $\mathscr{A}$ as the solution
$\xi(s)=(\mathbf{X}(s),\theta(s))$ to the system $\dfrac{\d }{\d
s}\xi(s)=\mathscr{A}(\xi(s)),$ i.e.
\begin{equation*}\dfrac{\d }{\d s}\mathbf{X}(s)=\ff(\mathbf{X}(s)), \:\:\:\dfrac{\d }{\d s}\theta(s)=1, \qquad (s \in \mathbb{R}),\end{equation*}
with $$\mathbf{X}(0)=\x,\quad \theta(0)=t.$$ It is clear that the
solution $\xi(s)$ to the above system is given by
$$\mathbf{X}(s)={\mathbf{\Phi}}(\x,s), \qquad \theta(s)=s+t,$$
and we  can define the flow of solution $\mathbf{\Psi}(\xi,s)=(
{\mathbf{\Phi}}(\x,s),s+t)$ associated to $\mathscr{A}$ and the
existence times of the characteristic curves of $Y$ are defined, for
any $\xi=(\x,t) \in \O_T$, as
$$\ell_{\pm}(\xi)=\inf\{s >0,\left({\mathbf{\Phi}}(\x,\pm s),\pm s+t\right) \notin
\O_T\}.$$ The flow $\mathbf{\Psi}(\cdot,\cdot)$ enjoys,
\textit{mutatis mutandis}, the properties listed in Proposition
\ref{Phiprop} and $\mu_T$ is invariant under $\mathbf{\Psi}$.
Moreover, since $\mathscr{A}$ is clearly Lipschitz continuous on
$\overline{\O_T}$, no characteristic of $Y$ can escape to infinity
in finite time. In other words, all characteristic curves of $Y$ now
have finite lengths. Indeed, if ${\mathbf{\Phi}}(\x,\pm s)$ does not
reach $\partial \O$, then the characteristic curve
$\mathbf{\Psi}(\xi,\pm s)$
 enters or leaves $\O_T$ through the bottom
$\O \times \{0\}$, or through the top $ \O \times \{T\}$ of it.
Precisely, it is easy to verify that for $\xi=(\x,t) \in\O_T$ we
have
$$\ell_{+}(\xi)=\tau_+(\x) \wedge (T-t) \quad \text{ and } \quad
\ell_-(\xi)=\tau_-(\x) \wedge t,$$ where $\wedge$ denotes minimum.
This clearly implies $\sup\{\ell_{\pm}(\xi)\;;\,\xi \in \O_T\,\}\leq T.$
Define now
$$\Sigma_{\pm,\,T}=\{\zeta \in \partial\O_T\,;\,\exists \xi \in\O_T\,\text{ such that }\:
\zeta=\mathbf{\Psi}(\xi,\pm \ell_{\pm}(\xi))\}.$$
 The definition of
$\Sigma_{\pm,\,T}$ is analogous to $\Gamma_{\pm}$ with the understanding
that now the charateristic curves correspond to the vector field
$\mathscr{A}$. In other words,
 $\Sigma_{-,\,T}$ (resp. $\Sigma_{+,\,T}$) is the subset of
$\partial\O_T$ consisting of all left (resp. right) limits of
characteristic curves of $\mathscr{A}$ in $\O_T$ whereas
$\Gamma_{-}$ (resp. $\Gamma_{+}$) is the subset of $\partial \O$
consisting of all left (resp. right) limits of characteristic curves
of $\ff$ in $\O.$ The main difference (and the interest of such a
lifting to $\O_T$) is the fact that \textit{each characteristic
curve of $\mathscr{A}$ does reach the boundaries $\Sigma_{\pm,\,T}$
in finite time}. The above formulae allow us to extend functions
$\ell_{\pm}$ to $\Sigma_{\pm,\,T}$ in the same way as we extended the
functions $\t_\pm$ to $\Gamma_\pm$. With these considerations, we
can represent, up to a set of zero measure, the phase space $\O_T$
as
\begin{equation}\begin{split}\label{parametre}
\O_{T}&=\{\mathbf{\Psi}(\xi,s)\,;\,\xi \in \Sigma_{-,\,T}\,,\,0 < s
<\ell_+(\xi)\}
\\
&=\{\mathbf{\Psi}(\xi,-s)\,;\,\xi \in \Sigma_{+,\,T}\,,\,0 < s
<\ell_-(\xi)\}
.\end{split}\end{equation} With this realization we can prove the
following:
\begin{propo} Let $T > 0$ be fixed. There are unique positive Borel measures $\d\nu_{\pm}$
on $\Sigma_{\pm,T}$ such that $\d\mu_T=\d\nu_+\otimes\d
s=\d\nu_-\otimes\d s.$
\end{propo}
\begin{proof} For any $\delta>0$, define $\mathscr{E}_\delta$ as the set of all bounded Borel subsets
$E$ of $\Sigma_{-,T}$ such that $\ell_+(\xi) > \delta$ for any $\xi
\in E$. Let us now \textit{fix} $E \in \mathscr{E}_\delta$. For all
$0 < \sigma \leq \delta$ put
$$E_{\sigma} = \{\mathbf{\Psi}(\xi,s)\,;\,\xi \in E ,0 < s \leq \sigma \}.$$ Clearly
$E_{\sigma}$ is a measurable subset of $\O_{T}$. Define the mapping
$h:\:\sigma \in (0,\delta] \mapsto h(\sigma)=\mu_T(E_{\sigma})$ with
$h(0) =0$. If $\sigma_1$ and $\sigma_2$ are two positive numbers
such that $ \sigma_1 + \sigma_2 \leq \delta$, then
$$E_{\sigma_1 + \sigma_2}\setminus E_{\sigma_1} = \{\mathbf{\Psi}(\xi,s)\,;\,\xi \in E ,\sigma_1 < s \leq \sigma_1 + \sigma_2
\}= \{\mathbf{\Psi}(\eta,\sigma_1)\,;\, \eta \in E_{\sigma_2}\}.$$
The properties of the flow $\mathbf{\Psi}$  (see Proposition
\ref{Phiprop}) ensure that the mapping $\eta \mapsto
\mathbf{\Psi}(\eta,\sigma_1)$ is one-to-one and measure preserving,
so that
$$\mu_T(E_{\sigma_1 + \sigma_2}\setminus E_{\sigma_1}) =
\mu_T(E_{\sigma_2})= h(\sigma_2).$$ Since $E_{\sigma_1 + \sigma_2} =
E_{\sigma_1} \cup (E_{\sigma_1 + \sigma_2}\setminus E_{\sigma_1}),$
we immediately obtain
\begin{equation}\label{cauchy}
h(\sigma_1 + \sigma_2) = h(\sigma_1) + h(\sigma_2) \qquad   \text{
for any } \qquad \sigma_1,\,\sigma_2 >0\,\text{ with }\sigma_1
+\sigma_2 \leq \delta .\end{equation} This is the well-known Cauchy
equation, though defined only on an interval of the real line. It
can be solved in a standard way using non-negativity instead of
continuity, yielding:
$$h(\sigma) = c_E \sigma \qquad \text{ for any } \quad 0 < \sigma
\leq \delta$$ where $c_E = h(\delta)/\delta$.
We define $\nu_-(E) = c_E$. It is not difficult to see that, with
the above procedure, the mapping $\nu_-(\cdot)$ defines a positive
measure on the ring $\mathscr{E}=\bigcup_{\delta
>0}\mathscr{E}_\delta$ of all the Borel subsets of $\Sigma_{-,T}$
on which the function $\ell_+(\xi)$ is bounded away from $0$. Such a
measure $\nu_-$ can be uniquely extended to the $\sigma$-algebra of
the Borel subsets of $\Sigma_{-,T}$ (see e.g. \cite[Theorem A, p.\
54]{Hal}). Consider now a Borel subset $E$ of $\Sigma_{-,T}$ and a
Borel subset $I$ of $\mathbb{R}^+$, such that for all $\xi \in E$
and $s \in I$ we have $0 < s < \ell_+(\xi)$. Then
$$E \times I = \{\mathbf{\Psi}(\xi,s)\,;\,\xi \in E , s \in I \} \subset
\O_{T}.$$ Thanks to the definition of $\nu_-(\cdot)$,
 we can state that $\mu_T(E \times I)= \nu_-(E)\mathrm{meas}(I)$
where meas$(I)$ denotes the linear  Lebesgue measure of $I \subset
\mathbb{R}$. This shows that $\d\mu_T = \d\nu_-\otimes\d s.$
Similarly we can define a measure $\nu_+$ on $\Sigma_{+,T}$ and
prove that $\d\mu_T=\d\nu_+\otimes\d s.$ The uniqueness of the
measures $\d\nu_{\pm}$ is then obvious.\end{proof}
\begin{nb} Note that the above construction of the Borel measures $\d\nu_{\pm}$ differs from that of \cite[Lemmas
XI.3.1 \& 3.2]{proto}, \cite[Propositions 7 \& 8]{beals} which ,
moreover, only apply when $\mu$ is absolutely continuous with
respect to the Lebesgue measure. Our construction is much more
general and can also be generalized to the case of a non--divergence
force field $\ff$, \cite{m2as}.\end{nb}
 Next, by  the cylindrical structure of $\O_T$, and the representation of $\Sigma_{\pm,T}$ as
$$\Sigma_{-,\,T}=\left(\Gamma_- \times (0,T) \right) \cup\O \times
\{0\}\quad \text{ and } \quad \Sigma_{+,\,T}=\left(\Gamma_+ \times
(0,T) \right) \cup\O \times \{T\},$$ the measures $\d\nu_{\pm}$ over
$\Gamma_{\pm} \times (0,T)$ can be written as $\d\nu_{\pm}=\d\mu_{\pm} \otimes\d
t$, where $\d\mu_{\pm}$ are Borel measures on $\Gamma_{\pm}$.  This leads to
the following
\begin{lemme}\label{defmu+-} There are unique positive Borel measures $\d\mu_{\pm}$
on $\Gamma_{\pm}$ such that, for any $f \in L^1(\O_T,\d\mu_T)$
\begin{equation}\label{beals1}\begin{split}
\int_{\O_T}f(\x,t)\d\mu_T(\x,t)&=\int_0^T\d
t\int_{\Gamma_+}\d\mu_+(\y)\int_0^{\t_-(\y) \wedge t}
f({\mathbf{\Phi}}(\y,-s),t-s)\d s\\
&\phantom{x}+\int_{\O}\d\mu(\x)\int_0^{\t_-(\x) \wedge
T}f(\mathbf{\Phi}(\x,-s),T-s)\d s,
\end{split}\end{equation} and
\begin{equation}\label{beals2}\begin{split}
\int_{\O_T}f(\x,t)\d\mu_T(\x,t)&=\int_0^T\d
t\int_{\Gamma_-}\d\mu_-(\y)\int_0^{\t_+(\y) \wedge (T-t)}
f({\mathbf{\Phi}}(\y,s),t+s)\d s\\
&\phantom{x}+\int_{\O}\d\mu(\x)\int_0^{\t_+(\x) \wedge
T}f(\mathbf{\Phi}(\x,s),s)\d s.
\end{split}\end{equation}\end{lemme}
The above fundamental result allows to compute integrals over the
cylindrical phase-space $\O_T$ through integration along the
characteristic curves. Let us now generalize it to the phase space
$\O$. Here the main difficulty stems from the fact that the
characteristic curves of the vector field $\ff$ are no longer
assumed to be of finite length. In order to extend Lemma
\ref{defmu+-} to possibly infinite existence times, first we prove
the following:
\begin{lemme}
Let $T > 0$ be fixed. Then, $\t_+(\x) < T$ for any $\x \in \O$ if
and only if $\t_-(\x) < T$ for any $\x \in \O.$\end{lemme}
\begin{proof} It is easy to see that $\t_+(\x) < T$ for any $\x \in \O$ is equivalent to $\t(\x) < T$ for any $\x \in \O$ and this is also equivalent to
$\t_-(\x) < T$ for  any $\x \in \O$.\end{proof}

 Hereafter, the support of a measurable function $f$
defined on $\O$ is defined as $\mathrm{Supp}f=\O \setminus
\mathbf{\omega}$  where $\mathbf{\omega}$ is the maximal open subset
of $\O$ on which $f$ vanishes $\d \mu$--almost everywhere.

\begin{propo}\label{Prop10.10} Let $f \in L^1(\O,\d \mu)$. Assume that there exists $\t_0>0$ such that
$\t_{\pm}(\x) < \t_0$ for any $\x \in \mathrm{Supp}(f).$ Then,
\begin{equation}\begin{split}
\int_{\O}f(\x)\d\mu(\x)&=\int_{\Gamma_+}\d\mu_+(\y)\int_0^{\t_-(\y)}f\left({\mathbf{\Phi}}(\y,-s)\right)\d
s\\
&=\int_{\Gamma_-}\d\mu_-(\y)\int_0^{\t_+(\y)}f({\mathbf{\Phi}}(\y,s))\d
s.
\end{split}\end{equation}
\end{propo}
\begin{proof} For any $T > \t_0$, define the  domain $\O_{T}=\O \times
(0,T)$. Since $T < \infty,$ it is clear that $f \in L^1(\O_T,\d\mu
\d t)$ and, by \eqref{beals1}, we get
\begin{multline*}
T\int_{\O}f(\x)\d\mu(\x)=\int_0^T\d t \int_{\Gamma_+}\d
\mu_+(\y)\int_0^{t \wedge \t_-(\y)}f({\mathbf{\Phi}}(\y,-s))\d s+\\
\int_{\O}\d\mu(\x)\int_0^{\t_-(\x)}f({\mathbf{\Phi}}(\x,-s))\d
s.\end{multline*} Since the formula is valid for any $T > \t_0$,
differentiating with respect to $T$ leads to the first assertion.
The second assertion is proved in the same way by using formula
\eqref{beals2}.
\end{proof}

To drop the finiteness assumption on $\t_{\pm}(\x)$, first we introduce
the  sets
$$\O_{\pm}=\{\x \in \O\,;\,\t_{\pm}(\x) < \infty\}, \qquad
\O_{\pm\infty}=\{\x \in \O\,;\,\t_{\pm}(\x) = \infty\},$$ and
$$\Gamma_{\pm\infty}=\{\y \in \Gamma_{\pm}\,;\, \t_{\mp}(\y) = \infty\}.$$
Then
\begin{propo}\label{prointegra} Given $f \in L^1(\O,\d \mu)$, one has \begin{equation}\label{10.47}
\int_{\O_{\pm}}f(\x)\d\mu(\x)
=\int_{\Gamma_\pm}\d\mu_\pm(\y)\int_0^{\t_\mp(\y)}f\left({\mathbf{\Phi}}(\y,\mp
s)\right)\d s,
\end{equation}
and
\begin{equation}\label{10.49}
\int_{\O_{\pm} \cap
\O_{\mp\infty}}f(\x)\d\mu(\x)=\int_{\Gamma_\pm\infty}\d\mu_\pm(\y)\int_0^{\infty}f\left({\mathbf{\Phi}}(\y,\mp
s)\right)\d s.
\end{equation}
\end{propo}
\begin{proof} Assume first $f \geq 0$. Let us fix $T >0$. It is clear that $\x \in \O$ satisfies
$\t_+(\x) < T$ if and only if $\x={\mathbf{\Phi}}(\y,-s),$ with $\y
\in \Gamma_+$ and $0 < s < T \wedge \t_-(\y)$. Then, by Proposition
\ref{Prop10.10},
$$\int_{\{\t_+(\x) <
T\}}f(\x)\d\mu(\x)=\int_{\Gamma_+}\d\mu_+(\y)\int_0^{T \wedge
\t_-(\y)} f({\mathbf{\Phi}}(\y,-s))\d s.$$ Since $f \geq 0$, the
inner integral is increasing with $T$ and, using the monotone
convergence theorem, we let $T \to \infty$ to get
$$\int_{\O_{+}}f(\x)\d\mu(\x)
=\int_{\Gamma_+}\d\mu_+(\y)\int_0^{\t_-(\y)}f\left({\mathbf{\Phi}}(\y,-s)\right)\d
s$$ which coincides with \eqref{10.47}.  We proceed in the same way
with integration on $\Gamma_-$ and get the second part of
\eqref{10.47}. Next we consider the set
$$\Delta=\{\x \in \O\,;\,\x={\mathbf{\Phi}}(\y,-s),\,\y \in
\O_{+\infty},\,0<s<T\}.$$ Proposition \ref{Prop10.10} asserts that
$$\int_{\Delta}f(\x)\d\mu(\x)=\int_{\O_{+\infty}}\d\mu_+(\y)\int_0^Tf({\mathbf{\Phi}}(\y,-s))\d
s.$$ Letting again $T \to \infty$, we get \eqref{10.49}. We extend
the results to arbitrary $f$ by linearity.\end{proof}

Finally, with the following, we show that it is possible to transfer
integrals over $\Gamma_-$ to $\Gamma_+$:
\begin{propo}
For any $\psi \in L^1(\Gamma_-,\d\mu_-)$,
\begin{equation}\label{10.51}
\int_{\Gamma_-\setminus
\Gamma_{-\infty}}\psi(\y)\d\mu_-(\y)=\int_{\Gamma_+\setminus
\Gamma_{+\infty}}\psi({\mathbf{\Phi}}(\mathbf{z},-\t_-(\mathbf{z})))\d\mu_+(\mathbf{z}).\end{equation}
\end{propo}
\begin{proof} For any $\epsilon > 0$, let $f_{\e}$ be the
function defined on $\O_+ \cap \O_-$ by
\begin{equation*}
\psi_{\e}(\x)=\begin{cases}
\dfrac{\psi({\mathbf{\Phi}}(\x,-\t_-(\x)))}{\t_+(\x)+\t_-(\x)}
\quad &\text{ if } \quad \t_-(\x)+\t_+(\x) > \e,\\
0 \quad &\text{ else}.\end{cases}
\end{equation*}
Since $\psi_{\e} \in L^1(\O_+ \cap \O_-,\d \mu)$, Eqs. \eqref{10.47}
and \eqref{10.49} give
\begin{equation*}\begin{split}
\int_{\O_+ \cap \O_-} \psi_{\e}(\x)\d\mu(\x)&=\int_{\{\t_+(\y) >
\e\}\setminus
\Gamma_{-\infty}}\d\mu_-(\y)\int_0^{\t_+(\y)}\psi(\y)\dfrac{\d
s}{\t_+(\y)}\\
&=\int_{\{\t_+(\y) > \e\}\setminus
\Gamma_{-\infty}}\psi(\y)\d\mu_-(\y).\end{split}\end{equation*} In
the same way, \begin{equation*}\begin{split} \int_{\O_+ \cap \O_-}
\psi_{\e}(\x)\d\mu(\x)&=\int_{\{\t_-(\y)
> \e\}\setminus
\Gamma_{+\infty}}\d\mu_+(\y)\int_0^{\t_-(\y)}\psi({\mathbf{\Phi}}(\y,-\t_-(\y)))\dfrac{\d
s}{\t_-(\y)}\\
&=\int_{\{\t_-(\y) > \e\}\setminus
\Gamma_{+\infty}}\psi({\mathbf{\Phi}}(\y,-\t_-(\y)))\d\mu_-(\y),\end{split}\end{equation*}
which leads to
$$\int_{\{\t_-(\y) > \e\}\setminus
\Gamma_{+\infty}}\psi({\mathbf{\Phi}}(\y,-\t_-(\y)))\d\mu_+(\y)=\int_{\{\t_+(\y)
> \e\}\setminus \Gamma_{-\infty}}\psi(\y)\d\mu_-(\y)$$
for any $\e > 0$. Passing to the limit as $\e \to 0$ we get the
conclusion.\end{proof}

We end this section with a technical result we shall need in the
sequel (see Lemma \ref{supp}):
\begin{propo}\label{support} Let $K$ be a compact subset of $\O$.
Denote
$$K_\pm:=\left\{\y \in \Gamma_\pm\;;\;\exists  t_0
\in \mathbb{R}\;\text{ such that }  {\mathbf{\Phi}}(\y,\pm t) \in
K\, \text{ for any } t \geq t_0\right\}.$$ Then $\mu_\pm(K_\pm)=0.$
\end{propo}
\begin{proof} Let $K$ be a fixed compact subset of $\O$. Applying
Eq. \eqref{10.47} or \eqref{10.49} to the function
$f(\x)=\chi_K(\x)$, one has \begin{equation}\label{muK} \infty >
\mu(K) \geq \int_{K_-}\d\mu_-(\y)\int_0^\infty \chi_K(
{\mathbf{\Phi}}(\y,t))\d t.\end{equation} By definition, if $\y \in
K_-$, then for some $t_0 \in \mathbb{R}$,
$\chi_K({\mathbf{\Phi}}(\y,t))=1$ for any $t \geq t_0$. Therefore,
$$\int_0^\infty \chi_K( {\mathbf{\Phi}}(\y,t)) = \infty, \qquad \forall
\y \in K_-.$$ Inequality \eqref{muK} implies   that $\mu_-(K_-)=0$.
One proves the result for $K_+$ in the same way.
\end{proof}
\section{The maximal transport operator and trace
results}\label{sec:maxi}  The results of the previous section allow
us to define the (maximal) transport operator $\T_{\mathrm{max}}$ as
the weak derivative along the characteristic curves.   To be
precise, let us define the space of \textit{test functions}
$\mathfrak{Y}$ as follows:

\begin{defi}[\textit{\textbf{Test--functions}}] Let $\mathfrak{Y}$ be the set of all measurable and bounded
functions $\psi :\O \to \mathbb{R}$ with compact support in $\O$ and such
that, for any $\x \in \O$, the mapping
$$s \in (-\t_-(\x),\t_+(\x)) \longmapsto \psi({\mathbf{\Phi}}(\x,s))$$
is continuously differentiable with
\begin{equation}\label{defi:deriva}\x \in \O \longmapsto \dfrac{\d}{\d
s}\psi({\mathbf{\Phi}}(\x,s))\bigg|_{s=0} \text{ measurable and
bounded}.\end{equation}
\end{defi}

\begin{nb}  Notice
that the class of test-functions $\mathfrak{Y}$ is not defined as a
subset of $L^\infty(\O,\d\mu)$; that is, we do not identify
functions equal $\mu$-almost everywhere. It is however a natural
question to know whether two test-functions coinciding $\mu$-almost
everywhere are such that there derivatives (defined by
\eqref{defi:deriva}) do coincide  $\mu$-almost everywhere. We
provide a positive answer to this question at the end of the paper
(see Appendix).
\end{nb}

An important property of test-functions is the following consequence of Proposition \ref{support}:
\begin{lemme}\label{supp} Let $\psi \in \mathfrak{Y}$ be given. For $\mu_\mp$-almost any $\y \in \Gamma_{\mp}$ there exists a
sequence $(t_n^\pm)_n$ (depending on $\y$) such that
$$\lim_{n \to \infty} t^\pm_n=\t_\pm(\y) \qquad \text{ and } \quad
\psi(\mathbf{\Phi}(\y,\pm t_n^\pm))=0 \quad \forall n \in \mathbb{N}.$$
\end{lemme}
\begin{proof} Let $\psi \in \mathfrak{Y}$ be given and let $K=\mathrm{Supp}(\psi)$. For any $\y \in \Gamma_-$
with $\t_+(\y) < \infty$ one has $\mathbf{\Phi}(\y,\t_+(\y)) \in
\Gamma_+$ and, since $K$ is compact in $\O$,
$\psi(\mathbf{\Phi}(\y,\t_+(\y))=0$ and the existence of a sequence
$(t_n^+)_n$ converging to $\t_+(\y)$ with the above property is
clear. Now, Proposition \ref{support} applied to $K$ shows that
there exists a set $\Gamma_-' \subset \Gamma_-$ with $\mu_-(\Gamma
\setminus \Gamma_-')=0$  and such that, for any $\y \in \Gamma_-'$,
there is a sequence $(t_n^+)_n$ converging to $\infty$ such that
$\mathbf{\Phi}(\y,t_n) \notin K$ for any $n \in \mathbb{N}$. This
proves the result. The statement for $\Gamma_+$ is proved in the
same way.\end{proof}

In the next step we define the transport operator
$(\T_\mathrm{max},\D(\T_\mathrm{max}))$.
\begin{defi}[\textit{\textbf{Transport operator $\T_\mathrm{max}$}}]
The domain of the maximal transport operator $\T_\mathrm{max}$ is
the set $\D(\T_\mathrm{max} )$ of all $f \in L^1(\O ,\d\mu)$ for
which there exists $g \in L^1(\O ,\d\mu)$ such that
$$\int_{\O }g(\x)\psi(\x)\d\mu(\x)=\int_{\O }f(\x)\dfrac{\d}{\d
  s}\psi({\mathbf{\Phi}}(\x,s))\bigg|_{s=0} \d\mu(\x),\qquad \qquad \forall \psi \in
\mathfrak{Y}.$$ In this case, $g=:\T_\mathrm{max} f.$
\end{defi}
\begin{nb}\label{nb:weakderiv} Of course, in
some weak sense, $\T_{\mathrm{max}}f=-\ff \cdot \nabla f$.
Precisely, for any $\varphi \in \Con^1_0(\O)$, the following formula
holds:
$$\int_{\O}\left(\ff(\x)\cdot \nabla
  \varphi(\x)\right)f(\x)\d\mu(\x)=
  \int_{\O}\T_{\mathrm{max}}f(\x)\varphi(\x)\d\mu(\x).$$
 \end{nb}

\subsection{\textbf{Fundamental representation formula: mild formulation}} Recall  that, if $f_1$ and $f_2$
are two functions defined over $\O$, we say that $f_2$ is a {\it
representative of} $f_1$ if $\mu\{\x \in \O\,;\,f_1(\x) \neq
f_2(\x)\}=0$, i.e. when $f_1(\x)=f_2(\x)$ for $\mu$-almost every $\x
\in \O $. The following fundamental result provides a characterization of  the
domain of $\D(\T_\mathrm{max})$:
\begin{theo}\label{representation}
Let $f \in L^1(\O,\mu)$. The following are equivalent:
\begin{enumerate}\item There exists $g \in L^1(\O,\mu)$ and
a representative $f^\sharp$ of $f$ such that, for $\mu$-almost every
$\x \in \O$ and any $-\tau_-(\x) < t_1 \leq t_2 < \tau_+(\x)$:
\begin{equation}\label{integralTM-}f^{\sharp}({\mathbf{\Phi}}(\x,t_1))-f^{\sharp}({\mathbf{\Phi}}(\x,t_2))
=\int_{t_1}^{t_2}
g({\mathbf{\Phi}}(\x,s))\d s.
\end{equation}
\item $f \in \D(\T_\mathrm{max})$. In this case, $g=\T_\mathrm{max} f$.
\end{enumerate}
\end{theo}

 The proof of the theorem is made of several steps. The difficult
 part of the proof is the  implication $(2) \implies (1)$. It is carried out through
several technical lemmas based
 upon \textit{\textbf{mollification along the characteristic curves}}  (recall
 that, whenever $\mu$ is not absolutely continuous with respect to
 the Lebesgue measure, no global convolution argument is available). Let us make precise what this is all about. Consider
a sequence $(\varrho_n)_n$  of one dimensional mollifiers supported
in $[0,1]$, i.e. for any $n \in \mathbb{N}$, $\varrho_n \in
\Con^{\infty}_0(\mathbb{R})$, $\varrho_n(s)=0$ if $s \notin
[0,1/n]$, $\varrho_n(s) \geq 0$  and $\int_0^{1/n}\varrho_n(s)\d
s=1.$ Then, for any $f \in L^1(\O,\d\mu)$, define the (extended)  mollification:
$$\varrho_n \diamond f(\x)=\int_0^{\t_-(\x)}
\varrho_n(s)f({\mathbf{\Phi}}(\x,-s))\d s.$$
As we shall see later, such a definition corresponds precisely to a
time convolution  over any characteristic curves (see e.g. \eqref{equamoll}).
Note that, with such a
definition, it is not clear {\it a priori} that $\varrho_n \diamond
f$ defines a measurable function, finite almost everywhere. It is
proved in the following that actually such a function is integrable.
\begin{lemme} Given $f \in L^1(\O,\d\mu)$, $\varrho_n \diamond f \in
  L^1(\O ,\d\mu)$ for any $n \in \mathbb{N}$. Moreover,
\begin{equation}\label{contrac}
\|\varrho_n \diamond f\| \leq \|f\|,\qquad \forall f \in
L^1(\O,\d\mu), n \in \mathbb{N}.\end{equation}
\end{lemme}
\begin{proof} One considers, for a given $f \in
L^1(\O,\d\mu)$,  the extension of $f$ by zero outside $\O$:
$$\overline{f}(\x)=f(\x),\qquad \forall \x \in \O,\qquad
\overline{f}(\x)=0 \quad \forall \x \in \mathbb{R}^N \setminus \O.$$
Then $\overline{f} \in L^1(\mathbb{R}^N,\d\mu).$ Let us consider the
 transformation:
$$\Upsilon\::\:(\x,s) \in \mathbb{R}^N \times \mathbb{R} \mapsto
\Upsilon(\x,s)=({\mathbf{\Phi}}(\x,-s),-s) \in \mathbb{R}^N \times
\mathbb{R}.$$ As a homeomorphism, $\Upsilon$ is measure preserving
for  pure Borel measures. It is also measure preserving for
completions of Borel measures (such as a Lebesgue measure) since
it is measure-preserving on Borel sets and the completion of a
measure is obtained by adding to the Borel $\sigma$-algebra all
sets contained in a measure-zero Borel sets, see \cite[Theorem
13.B, p.~55]{Hal}. Then, according to \cite[Theorem 39.B,
p.~162]{Hal}, the mapping
$$(\x,s) \in \mathbb{R}^N \times \mathbb{R} \mapsto
\overline{f}({\mathbf{\Phi}}(\x,-s))$$ is measurable as the
composition of $\Upsilon$ with the measurable function $(\x,s)
\mapsto \overline{f}(\x)$. Define now $\Lambda=\{(\x,s)\,;\,\x \in
\O,\,0 < s < \t_-(\x)\},$ $\Lambda$ is a measurable subset of
$\mathbb{R}^N \times \mathbb{R}$. Therefore, the mapping
$$(\x,s) \in \mathbb{R}^N \times \mathbb{R} \longmapsto
\overline{f}({\mathbf{\Phi}}(\x,-s))\chi_{\Lambda}(\x,s)
\varrho_n(s)$$ is measurable. Since $\varrho_n$ is compactly
supported, it is also integrable over $\mathbb{R}^N \times
\mathbb{R}$ and, according to Fubini's Theorem
\begin{equation*}
[\varrho_n \diamond f](\x):= \ds
\int_{\mathbb{R}}\overline{f}({\mathbf{\Phi}}(\x,-s))\chi_{\Lambda}(\x,s)
\varrho_n(s)\d s = \ds
\int_0^{\tau_-(\x)}\varrho_n(s)f({\mathbf{\Phi}}(\x,-s))\d
s\end{equation*} is finite for almost every $\x \in \O$ the
 and the associated application $\varrho_n \diamond f$ is
 integrable.

 Let us prove now \eqref{contrac}. Since $|\varrho_n \diamond f| \leq \varrho_n \diamond |f|$, to
  show that $\varrho_n \diamond f \in L^1(\O,\d\mu)$, it suffices to deal
  with a \textit{nonnegative function} $f \in L^1(\O,\d\mu)$. One sees easily that,
  for any $\y \in \Gamma_-$ and any $0 < t < \t_+(\y)$,
\begin{equation}\label{equamoll}(\varrho_n \diamond f)({\mathbf{\Phi}}(\y,t))=\int_0^t
\varrho_n(s)f({\mathbf{\Phi}}(\y,t-s))\d s=\int_0^t \varrho_n(t-s)
f({\mathbf{\Phi}}(\y,s))\d s.\end{equation} Thus,
\begin{equation*}\begin{split}
\int_0^{\t_+(\y)}[\varrho_n \diamond f]({\mathbf{\Phi}}(\y,t))\d t
&=\int_0^{\t_+(\y)}\d s
\int_s^{\t_+(\y)}\varrho_n(s)f({\mathbf{\Phi}}(\y,t-s))\d t\\
&=\int_0^{\t_+(\y) \wedge 1/n}\varrho_n(s)\d
s\int_0^{\t_+(\y)-s}f({\mathbf{\Phi}}(\y,r))\d r.
\end{split}\end{equation*}
Therefore, \begin{equation*}\begin{split} 0 \leq
\int_0^{\t_+(\y)}[\varrho_n \diamond f]({\mathbf{\Phi}}(\y,t))\d t
&\leq \int_0^{1/n}\varrho_n(s)\d
s\int_0^{\t_+(\y)}f({\mathbf{\Phi}}(\y,r))\d
r\\
&=\int_0^{\t_+(\y)}f({\mathbf{\Phi}}(\y,r))\d r,\qquad \forall \y
\in \Gamma_-,\:n \in \mathbb{N}\end{split}\end{equation*} so that
$$\int_{\Gamma_-}\d\mu_-(\y)\int_0^{\t_+(\y)}[\varrho_n \diamond
f]({\mathbf{\Phi}}(\y,t))\d t \leq
\int_{\Gamma_-}\d\mu_-(\y)\int_0^{\t_+(\y)}
f({\mathbf{\Phi}}(\y,r))\d r.$$ This proves, thanks to Proposition
\ref{prointegra}, that
\begin{equation}\label{intO_-}\int_{\O_-}[\varrho_n \diamond f]\d\mu \leq
\int_{\O_-}f\d\mu.\end{equation} Now, in the same way:
\begin{equation*}\begin{split}
\int_{\O_+ \cap \O_{-\infty}}&[\varrho_n \diamond f](\x)\d\mu(\x)=
\int_{\Gamma_{+\infty}}\d\mu_+(\y)\int_0^\infty [\varrho_n \diamond
f]({\mathbf{\Phi}}(\y,-t))\d t\\
&=\int_{\Gamma_{+\infty}}\d\mu_+(\y)\int_0^\infty\d
t\int_0^\infty\varrho_n(s)f({\mathbf{\Phi}}(\y,-s-t))\d s\\
&=\int_{\Gamma_{+\infty}}\d\mu_+(\y)\int_0^\infty\d
t\int_t^\infty\varrho_n(r-t)f({\mathbf{\Phi}}(\y,-r))\d
r.\end{split}\end{equation*} so that
\begin{equation*}\begin{split}
\int_{\O_+ \cap \O_{-\infty}}[\varrho_n \diamond f](\x)\d\mu(\x)
&=\int_{\Gamma_{+\infty}}\d\mu_+(\y)\int_0^\infty
f({\mathbf{\Phi}}(\y,-r))\d
r\int_0^r\varrho_n(r-t)\d t\\
&\leq \int_{\Gamma_{+\infty}}\d\mu_+(\y)\int_0^\infty
f({\mathbf{\Phi}}(\y,-r))\d r
\end{split}\end{equation*} i.e.
\begin{equation}\label{intO_+}
\int_{\O_+ \cap \O_{-\infty}}\varrho_n \diamond f(\x)\d\mu(\x)\leq
\int_{\O_+ \cap \O_{-\infty}} f(\x)\d\mu(\x).\end{equation} Finally
\begin{equation*}\begin{split}
\int_{\O_{+\infty} \cap \O_{-\infty}}[\varrho_n \diamond
f](\x)\d\mu(\x)&= \int_{\O_{+\infty} \cap
\O_{-\infty}}\d\mu(\x)\int_0^\infty\varrho_n(s)f({\mathbf{\Phi}}(\x,-s))\d s\\
&=\int_0^\infty\varrho_n(s)\d s\int_{\O_{+\infty}\cap
\O_{-\infty}}f({\mathbf{\Phi}}(\x,-s))\d \mu(\x).
\end{split}\end{equation*} Now, from Assumption \ref{ass:h2}, for any $s
\geq 0,$
$$\int_{\O_{+\infty}\cap
\O_{-\infty}}f({\mathbf{\Phi}}(\x,-s))\d
\mu(\x)=\int_{\O_{+\infty}\cap \O_{-\infty}}f(\x)\d \mu(\x),$$ so
that
\begin{equation}\label{intOinfty} \int_{\O_{+\infty} \cap
\O_{-\infty}}[\varrho_n \diamond
f](\x)\d\mu(\x)=\int_{\O_{+\infty}\cap \O_{-\infty}}f(\x)\d
\mu(\x).\end{equation} Combining \eqref{intO_-}, \eqref{intO_+} and
\eqref{intOinfty}, one finally gets $\|\varrho_n \diamond f\|  \leq
\|f\| .$\end{proof}

As it is the case for classical convolution, the family $(\varrho_n
\diamond f)_n$ approximates $f$ in $L^1$-norm:
\begin{propo}\label{approx} Given $f \in L^1(\O,\d\mu)$,
\begin{equation}\label{convergconv}
\lim_{n \to \infty}\int_{\O }\bigg|(\varrho_n \diamond f)(\x)
-f(\x)\bigg|\d\mu(\x)=0.\end{equation}
\end{propo}
\begin{proof}
 According to \eqref{contrac} and from the density of $\mathscr{C}_0(\O)$ in $L^1(\O,\d\mu)$, it
suffices to prove the result for any $f$ continuous over $\O$ and
compactly supported. Splitting $f$ into positive  and negative
parts, $f=f^+-f^-$, one can also assume $f$ to be nonnegative.  From
the continuity of both $f$ and $\mathbf{\Phi}(\cdot,\cdot)$, one has
$$\mathscr{K}_n:=\mathrm{Supp}(\varrho_n \diamond f)=\overline{\bigg\{\x \in \O\,,\,\exists s_0
\in \mathrm{Supp}(\varrho_n) \text{ such that } \mathbf{\Phi}(\x,-s_0) \in
\mathrm{Supp}(f)\bigg\}}.$$  Moreover, it is easily seen that
$\mathscr{K}_{n+1} \subset \mathscr{K}_n$ for any $n \geq 1$.
Finally, it is clear that $$\mathscr{K}_1 \subset \{\x \in
\overline{\O}\,;\exists \y \in \mathrm{Supp}(f)\,\text{ with } |\x-\y|
\leq d\}$$ where $d=\sup\{|\mathbf{\Phi}(\x, s)-\x|\,;\,0 \leq s
\leq 1\,; \x \in \mathrm{Supp}(f)\} < \infty.$ Therefore,
$\mathscr{K}_1$ is compact. Set now
$$\mathcal{O}_n:=\mathscr{K}_n \cup \mathrm{Supp}(f) \qquad \text{
and } \qquad
 \mathcal{O}_n^-=\{\x \in \mathcal{O}_n\,;\,\t_-(\x) <1/n\}.$$
Noticing that $\mu(\mathcal{O}_1)$ is finite, one can see easily
that $\lim_n\mu(\mathcal{O}_n^-)=0$. Since  $\sup_{\x \in \O}
|\varrho_n \diamond f(\x)| \leq \sup_{\x \in \O} |f(\x)|,$ for any
$\varepsilon > 0$, there exists $n_0 \geq 1$ such that
$$\int_{\mathcal{O}_n^-}|f(\x)|\d\mu(\x) \leq \varepsilon, \quad \text{ and }
\quad \int_{\mathcal{O}_n^-}|\varrho_n \diamond f(\x)|\d\mu(\x) \leq
\varepsilon \qquad \forall n \geq n_0.$$ Now, noticing that $
\mathrm{Supp}(\varrho_n \diamond f -f) \subset \mathcal{O}_n $, one
has for any $n \geq n_0$,
$$\int_{\O}|\varrho_n \diamond f-f|\d\mu=\int_{\mathcal{O}_n}|\varrho_n
\diamond f-f| \leq 2\varepsilon + \int_{\mathcal{O}_n \setminus
  \mathcal{O}_n^-}|\varrho_n \diamond f-f|\d\mu.$$
For any $\x \in \mathcal{O}_n \setminus \mathcal{O}_n^-$, since
$\varrho$ is supported in $[0,1/n]$, one has
\begin{equation*}\begin{split}
[\varrho_n \diamond f](\x)-f(\x) &=\int_0^{1/n}\varrho_n(s)
 f({\mathbf{\Phi}}(\x, -s))\d
s-f(\x)\\
&=\int_0^{1/n}\varrho_n(s)\left(f({\mathbf{\Phi}}(\x,
-s))-f(\x)\right)\d s.\end{split}\end{equation*} Note that, thanks
to Gronwall's lemma,
$$|\mathbf{\Phi}(\x,-s)-\x| \leq \frac{L}{\kappa}(\exp(k s)-1) \leq
\frac{L}{\kappa}(\exp(\kappa/n)-1),\qquad \forall \x \in
\mathcal{O}_1, s \in (0,1/n)$$ where $L=\sup\{|\ff(\x)|,\x \in
\mathcal{O}_1\}.$ Since $f$ is   uniformly continuous on
$\mathcal{O}_1$, it follows that
$$\lim_{n \to \infty} \sup \bigg\{|f(\mathbf{\Phi}(\x,-s)-f(\x)|\,;\,\x \in
\mathcal{O}_1,\;s \in (0,1/n)\bigg\}=0$$ from which we deduce  that
there exists some $n_1 \geq 0$, such that $|\varrho_n \diamond
f(\x)-f(\x)| \leq \varepsilon$ for any $\x \in \mathcal{O}_n
\setminus \mathcal{O}_n^-$ and any $n \geq n_1$. One obtains then,
for any $n \geq n_1$,
$$\int_{\O}|\varrho_n \diamond f-f|\d\mu \leq 2\varepsilon +\varepsilon
\mu(\mathcal{O}_n \setminus \mathcal{O}_n^-) \leq 2\varepsilon
+\varepsilon \mu(\mathcal{O}_1)$$ which proves the
result.\end{proof}

We saw that, for a given $f \in L^1(\O,\d\mu)$,
$\varrho_n \diamond f$ is also integrable ($n \in \mathbb{N}$).
Actually, we shall see that $\varrho_n \diamond f$ is even more
regular than $f$:
\begin{lemme}\label{moll} Given $f \in L^1(\O,\d\mu)$, set $f_n=\varrho_n \diamond f$, $n \in \mathbb{N}$.
Then, $f_n \in \D(\T_\mathrm{max})$ with
\begin{equation*}
[\T_\mathrm{max}
f_n](\x)=-\int_0^{\t_-(\x)}\varrho_n'(s)f({\mathbf{\Phi}}(\x,-s))\d
s, \qquad \x \in \O.
\end{equation*}
\end{lemme}
\begin{proof} Set $g_n(\x)=-\int_0^{\t_-(\x)}\varrho_n'(s)f({\mathbf{\Phi}}(\x,-s))\d
s,$ $\x \in \O.$ It is easy to see that $g_n \in L^1(\O,\d\mu)$.
Now, given $\psi \in \mathfrak{Y}$, let us consider the quantity
$$I=\int_{\O }f_n(\x)\dfrac{\d}{\d
  s}\psi({\mathbf{\Phi}}(\x,s))\bigg|_{s=0} \d\mu(\x).$$
One has to prove that $I=\int_{\O}g_n(\x)\psi(\x)\d\mu(\x)$. We split
the above integral over $\O$ into three integrals $I_-$, $I_+$ and
$I_\infty$ over $\O_-$, $\O_{-\infty}\cap \O_+$ and
$\O_{+\infty}\cap \O_{-\infty}$ respectively.
 Recall that, for any $\x
\in \O_-$, there is some $\y \in \Gamma_-$ and some $t \in
(0,\t_+(\y))$ such that $\x={\mathbf{\Phi}}(\y,t)$. In such a case
\begin{equation}\label{lip}
\dfrac{\d}{\d
s}\psi({\mathbf{\Phi}}(\x,s))\bigg|_{s=0}=\dfrac{\d}{\d
t}\psi({\mathbf{\Phi}}(\y,t)).\end{equation}  Then, according to
Prop. \ref{prointegra} and Eq. \eqref{equamoll}:
\begin{equation}\label{i-}\begin{split}
I_-&=\int_{\Gamma_-}\d\mu_-(\y)\int_0^{\t_+(\y)}f_n({\mathbf{\Phi}}(\y,t))\frac{\d}{\d
  t}\psi({\mathbf{\Phi}}(\y,t))\d t\\
  &=\int_{\Gamma_-}\d\mu_-(\y)\int_0^{\t_+(\y)} \frac{\d}{\d
  t}\psi({\mathbf{\Phi}}(\y,t))\d
  t\int_0^t\varrho_n(t-s)f({\mathbf{\Phi}}(\y,s))\d s\\
  &=\int_{\Gamma_-}\d\mu_-(\y)\int_0^{\t_+(\y)} f({\mathbf{\Phi}}(\y,s))\d s \int_s^{\t_+(\y)}\frac{\d}{\d
  t}\psi({\mathbf{\Phi}}(\y,t))
   \varrho_n(t-s) \d t.
  \end{split}\end{equation}
 Let us now investigate more carefully this last integral. Let $\y \in \Gamma_-$ be fixed.
 If $\t_+(\y) < \infty$ then, since $\psi$ is compactly supported, we have $\psi(\mathbf{\Phi}(\y,\t_+(\y)))=0$ and integration by part (together with $\varrho_n(0)=0$) leads to
 $$\int_s^{\t_+(\y)}\frac{\d}{\d
  t}\psi({\mathbf{\Phi}}(\y,t))
   \varrho_n(t-s) \d t=-\int_s^{\t_+(\y)}\varrho_n'(t-s)\psi(\mathbf{\Phi}(\y,t))\d t.$$
If now $\t_+(\y) > \infty$, then, since $\varrho_n$ is supported in $[0,1/n]$, one has
\begin{equation*}\begin{split}
      \int_s^{\t_+(\y)}\frac{\d}{\d
  t}\psi({\mathbf{\Phi}}(\y,t))
   \varrho_n(t-s) \d t&=\int_s^{s+1/n}\frac{\d}{\d
  t}\psi({\mathbf{\Phi}}(\y,t))
   \varrho_n(t-s) \d t\\
   &=-\int_s^{\t_+(y)}\varrho_n'(t-s)\psi(\mathbf{\Phi}(\y,t))\d t
                 \end{split}\end{equation*}
Finally, we obtain,
\begin{equation*}\begin{split}I_-&=-\int_{\Gamma_-}\d\mu_-(\y)\int_0^{\t_+(\y)} f({\mathbf{\Phi}}(\y,s))\d s  \int_s^{\t_+(\y)} \psi({\mathbf{\Phi}}(\y,t))
   \varrho_n'(t-s) \d t \\
&=-\int_{\Gamma_-}\d\mu_-(\y)\int_0^{\t_+(\y)}
\psi({\mathbf{\Phi}}(\y,t))\d
  t\int_0^t\varrho_n'(s)f({\mathbf{\Phi}}(\y,t- s))\d s.
   \end{split}\end{equation*}
Using again Prop. \ref{prointegra}, we finally get
$$I_-=\int_{\O_-}g_n(\x)\psi(\x)\d\mu(\x).$$
One proves in the same way that
$$I_+=\int_{\O_+ \cap \O_{-\infty} }f_n(\x)\dfrac{\d}{\d
  s}\psi({\mathbf{\Phi}}(\x,s))\bigg|_{s=0} \d\mu(\x)=\int_{\O_+\cap
  \O_{-\infty}}g_n(\x)\psi(\x)\d\mu(\x).$$
It remains to consider $I_\infty=\int_{\O_{+\infty} \cap
\O_{-\infty} }f_n(\x)\frac{\d}{\d
  s}\psi({\mathbf{\Phi}}(\x,s))\big|_{s=0} \d\mu(\x).$ One has
\begin{equation*}\begin{split}
I_\infty&=\int_{\O_{+\infty} \cap \O_{-\infty}}\frac{\d}{\d
  s}\psi({\mathbf{\Phi}}(\x,s))\bigg|_{s=0}\d\mu(\x)\int_0^\infty
  \varrho_n(t)f({\mathbf{\Phi}}(\x,-t))\d t\\
  &=\int_0^\infty \varrho_n(t)\d t\int_{\O_{+\infty} \cap \O_{-\infty}}\frac{\d}{\d
  s}\psi({\mathbf{\Phi}}(\x,s))\bigg|_{s=0}f({\mathbf{\Phi}}(\x,-t))\d
  \mu(\x).
\end{split}\end{equation*}
For any $\x \in \O_{+\infty} \cap \O_{-\infty}$ and any $t \geq 0$,
setting $\y={\mathbf{\Phi}}(\x,-t)$, one has $\y \in \O_{-\infty}
\cap \O_{+\infty}$ and $\frac{\d}{\d
  s}\psi({\mathbf{\Phi}}(\x,s))\big |_{s=0}=\frac{\d}{\d
  t}\psi({\mathbf{\Phi}}(\y,t))$ from which Liouville's Theorem (Assumption \ref{ass:h2}) yields
$$\int_{\O_{+\infty} \cap \O_{-\infty}}\frac{\d}{\d
  s}\psi({\mathbf{\Phi}}(\x,s))\bigg|_{s=0}f({\mathbf{\Phi}}(\x,-t))\d
  \mu(\x)=\int_{\O_{+\infty} \cap \O_{-\infty}}\frac{\d}{\d
  t}\psi({\mathbf{\Phi}}(\y,t)) f(\y)\d
  \mu(\y).$$
Therefore,
\begin{eqnarray*}I_\infty&=&\int_{\O_{+\infty} \cap \O_{-\infty}} f(\y)\d
  \mu(\y)\int_0^\infty \varrho_n(t)\frac{\d}{\d
  t}\psi({\mathbf{\Phi}}(\y,t))\d t\\&=&-\int_{\O_{+\infty} \cap \O_{-\infty}} f(\y)\d
  \mu(\y)\int_0^\infty \varrho_n'(t) \psi({\mathbf{\Phi}}(\y,t))\d t\\
  &=&-\int_0^\infty \varrho_n'(t)\d t\int_{\O_{+\infty} \cap \O_{-\infty}} f(\y)\psi({\mathbf{\Phi}}(\y,t))\d
  \mu(\y).
\end{eqnarray*} Arguing as above, one can "turn back" to the $\x$
variable  to get
$$\int_{\O_{+\infty} \cap \O_{-\infty}} f(\y)
\psi({\mathbf{\Phi}}(\y,t))\d\mu(\y)=\int_{\O_{+\infty} \cap
\O_{-\infty}} f({\mathbf{\Phi}}(\x,-t))\psi(\x)\d\mu(\x),$$ i.e.
$$I_\infty=-\int_{\O_{+\infty} \cap \O_{-\infty}}\psi(\x)\d\mu(\x)\int_0^\infty \varrho_n'(t)f({\mathbf{\Phi}}(\x,-t))\d
t=\int_{\O_{+\infty} \cap \O_{-\infty}}\psi(\x)g_n(\x)\d\mu(\x)$$
and the Lemma is proven.\end{proof}
\begin{nb}
Notice that Proposition \ref{approx} together with Lemma \ref{moll}
prove that $\D(\T_\mathrm{max})$ is a dense subset of
$L^1(\O,\d\mu)$.
\end{nb} Now, whenever $f \in \D(\T_\mathrm{max} )$, one has the following
more precise result:
\begin{propo}\label{moll2}
If $f \in \D(\T_\mathrm{max}  )$, then
\begin{equation}\label{claimdiamond}
[\T_\mathrm{max}  (\varrho_n \diamond f)](\x)=[\varrho_n \diamond
\T_\mathrm{max} f](\x),\qquad \qquad (\x \in \O\,,\,n \in
\mathbb{N}).\end{equation}\end{propo}

Before proving this result, we need the following very simple lemma:
\begin{lemme}\label{lemme:prop4.5} For any $\psi \in \mathfrak{Y}$ and any $n \in
\mathbb{N}$, define
$$\chi_n(\x)=
\int_0^{\t_+(\x)}\varrho_n(s)\psi({\mathbf{\Phi}}(\x,s))\d s, \qquad
\x \in \O.$$ Then, $\chi_n$ belongs to $\mathfrak{Y}$.
\end{lemme}
\begin{proof} Since $\t_+$ is measurable and $\varrho_n$ is
compactly supported, it is easy to see that $\chi_n$ is measurable
and bounded over $\O$. Now, for any $\x \in \O$, and any $t \in
(\t_-(\x),\t_+(\x))$, one has
$$\chi_n({\mathbf{\Phi}}(\x,t))=\int_t^{\t_+(\x)}
\varrho_n(s-t)\psi({\mathbf{\Phi}}(\x,s))\d s.$$
It is clear then from the properties of $\varrho_n$ that the mapping
$t \in (-\t_-(\x),\t_+(\x)) \mapsto \chi_n({\mathbf{\Phi}}(\x,t))$
is continuously differentiable  with
\begin{equation}\label{derivativechin}
\dfrac{\d }{\d t}
\chi_n({\mathbf{\Phi}}(\x,t))=-\int_t^{\t_+(\x)}\varrho_n'(s-t)\psi({\mathbf{\Phi}}(\x,s))\d
s =\int_t^{\t_+(\x)} \varrho_n(s-t) \dfrac{\d }{\d
s}\left[\psi({\mathbf{\Phi}}(\x,s))\right]\d s.
\end{equation}
In particular, for $t=0$, one gets
$$\dfrac{\d }{\d t}
\chi_n({\mathbf{\Phi}}(\x,t))\bigg|_{t=0}=-\int_0^{\t_+(\x)}\varrho_n'(s)\psi({\mathbf{\Phi}}(\x,s))\d
s.$$ Since $\varrho_n'$ is compactly supported and $\psi \in
\mathfrak{Y}$, the application $\x\in \O \longmapsto \dfrac{\d }{\d
t} \chi_n({\mathbf{\Phi}}(\x,t))\big|_{t=0}$ is measurable and
bounded.\end{proof}

\begin{proof1} We use the notations of Lemma \ref{moll}. Since $\varrho_n
\diamond \T_\mathrm{max}  f \in L^1(\O,\d\mu)$, it suffices to show
that
$$\int_{\O  }f_n(\x)\dfrac{\d}{\d
s}\psi({\mathbf{\Phi}}(\x,s))\bigg|_{s=0}\d\mu(\x)=\int_{\O}\psi(\x)
[\varrho_n \diamond \T_\mathrm{max}  f](\x)\d\mu(\x), \qquad \forall
\psi \in \mathfrak{Y}.$$ Here again, we shall deal separately with
the integrals over $\O_-$, $\O_{+}\cap \O_{-\infty}$ and
$\O_{+\infty} \cap \O_{-\infty}.$
 Let $\chi_n$ be defined as in Lemma \ref{lemme:prop4.5}, as we
 already saw it (see \eqref{derivativechin}), for any $\y \in \Gamma_-,$ and any $0<
s<\t_+(\y)$,
$\frac{\d}{\d
s}\chi_n({\mathbf{\Phi}}(\y,s))=\int_s^{\t_+(\y)}\varrho_n(t-s)\frac{\d}{\d
t}[\psi({\mathbf{\Phi}}(\y,t))]\d t.$ Consequently, according to
\eqref{i-},
\begin{equation*}\begin{split} \int_{\O_-
} f_n(\x)&\dfrac{\d}{\d
s}\psi({\mathbf{\Phi}}(\x,s))\bigg|_{s=0}\d\mu(\x)=\int_{\Gamma_-}\d\mu(\y)\int_0^{\t_+(\y)}f({\mathbf{\Phi}}(\y,r))
\dfrac {\d}{\d r}\chi_n({\mathbf{\Phi}}(\y,r))\d r\\
&=\int_{\O_-}f(\x)\dfrac{\d}{\d
s}\chi_n({\mathbf{\Phi}}(\x,s))\bigg|_{s=0}\d\mu(\x)
=\int_{{\O}_-}\chi_n(\x)[\T_\mathrm{max} f](\x)\d\mu(\x)
\end{split}\end{equation*} where, for the two last identities, we
used \eqref{lip} and the fact that $\chi_n \in \mathfrak{Y}.$ Now,
using Prop. \ref{prointegra}
\begin{equation*} \begin{split}\int_{\O_-}&\chi_n(\x)[\T_\mathrm{max}  f](\x)\d\mu(\x)=\int_{\O_-}[\T_\mathrm{max}  f](\x)\d\mu(\x)
\int_0^{\t_+(\x)}\varrho_n(r)\psi({\mathbf{\Phi}}(\x,r))\d r\\
&=\int_{\Gamma_-}\d\mu_-(\y)\int_0^{\t_+(\y)}\psi({\mathbf{\Phi}}(\y,s))\d
s\int_0^s\varrho_n(s-t)[\T_\mathrm{max} f]({\mathbf{\Phi}}(\y,t))\d
t.
\end{split}\end{equation*}
Therefore, Eq. \eqref{equamoll} leads to
\begin{equation*}\label{TMO-}\begin{split}
\int_{\O_-} \chi_n(\x)[\T_\mathrm{max}
f](\x)\d\mu(\x)&=\int_{\Gamma_-}\d\mu_-(\y)\int_0^{\t_+(\y)}\psi({\mathbf{\Phi}}(\y,s))[\varrho_n
\diamond \T_\mathrm{max} f]({\mathbf{\Phi}}(\y,s))\d s\\
&=\int_{\O_-}\psi(\x)\left[\varrho_n \diamond
\T_{\mathrm{max}}f\right](\x)\d\mu(\x).
\end{split}\end{equation*} The integrals over $\O_+
\cap \O_{-\infty}$ and $\O_{-\infty}\cap \O_{+\infty}$ are evaluated
in the same way.
\end{proof1}

We are in position to prove the following
\begin{propo} Let $f \in L^1(\O,\d\mu)$ and $f_n=\varrho_n \diamond f$, $n \in \mathbb{N}$. Then, for $\mu_-$-- a. e.
$\y \in \Gamma_-$,
\begin{equation}\label{fny-}
f_n({\mathbf{\Phi}}(\y,s))-f_n({\mathbf{\Phi}}(\y,t))=\int_s^t
[\T_\mathrm{max} f_n]({\mathbf{\Phi}}(\y,r))\d r \qquad \forall 0 <
s < t < \t_+(\y).\end{equation} In the same way, for almost every
$\z \in \Gamma_+$,
$$f_n({\mathbf{\Phi}}(\z,-s))-f_n({\mathbf{\Phi}}(\z,-t))=\int_s^t \T_\mathrm{max}
f_n({\mathbf{\Phi}}(\z,-r))\d r, \qquad \forall 0 < s < t <
\t_-(\z).$$ \label{prop312}
\end{propo}
\begin{proof} We focus only on \eqref{fny-}, the second assertion following the same lines. Since $f \in L^1(\O_-,\d\mu)$,
  Proposition  \ref{prointegra} implies that the integral
$\int_0^{\t_+(\y)} |f({\mathbf{\Phi}}(\y,r))|\d r$ exists and is
finite for $\mu_-$-almost every $\y \in \Gamma_-$. Therefore, for
$\mu_-$-almost every $\y \in \Gamma_-$ and  any $0 < t < \t_+(\y)$,
the quantities $\int_0^t \varrho_n(t-s) f({\mathbf{\Phi}}(\y,s))\d
s$ and $\int_0^t \varrho_n'(t-s) f({\mathbf{\Phi}}(\y,s))\d s$ are
well-defined and finite. Moreover, thanks to Eq. \eqref{equamoll}
Lemma \ref{moll}, they are respectively equal to
$f_n({\mathbf{\Phi}}(\y,t))$ and $[\T_\mathrm{max}
f_n]({\mathbf{\Phi}}(\y,t)).$ In particular, the mapping $t \in
(0,\t_+(\y)) \mapsto [\T_\mathrm{max} f_n]({\mathbf{\Phi}}(\y,t))
\in \mathbb{R}$ is continuous. It is then easy to see that, for any
$0 < s < t < \t_+(\y)$
\begin{equation*}\begin{split}
\int_s^t &[\T_\mathrm{max} f_n]( {\mathbf{\Phi}}(\y,r))\d
r=-\int_s^t \d r \int_0^r
\varrho_n'(r-u) f({\mathbf{\Phi}}(\y,u))\d  u\\
&=-\int_0^s f({\mathbf{\Phi}}(\y,u))\d u \int_s^t \varrho_n'(r-u) \d
r - \int_s^t f({\mathbf{\Phi}}(\y,u))\d u \int_u^t \varrho_n'(r-u)
\d
r\\
&=-\int_0^t f({\mathbf{\Phi}}(\y,u))\varrho_n(t-u)\d u+\int_0^s
f({\mathbf{\Phi}}(\y,u))\varrho_n(s-u)\d u,
\end{split}
\end{equation*}
which is nothing but \eqref{fny-}.
\end{proof}
As a consequence, one gets the following result :
\begin{propo}\label{f+-}
For any $f \in \D(\T_\mathrm{max})$, there exists some functions
$\widetilde{f}_{\pm} \in L^1(\O_\pm,\d\mu)$ such that
$\widetilde{f}_\pm(\x)=f(\x)$ for $\mu$- almost every $\x \in
\O_\pm$ and, for $\mu_-$--almost every $\y \in \Gamma_-$:
\begin{equation}\label{fy-}
\widetilde{f }_-({\mathbf{\Phi}}(\y,s))-\widetilde{f
}_-({\mathbf{\Phi}}(\y,t))=\int_s^t [\T_\mathrm{max} f]
({\mathbf{\Phi}}(\y,r))\d r \qquad \forall 0 < s < t <
\t_+(\y),\end{equation} while, for $\mu_+$--almost every $\z \in
\Gamma_+$:
\begin{equation*}
\widetilde{f }_+({\mathbf{\Phi}}(\z,-s))-\widetilde{f
}_+({\mathbf{\Phi}}(\z,-t))=\int_s^t [\T_\mathrm{max} f ]
({\mathbf{\Phi}}(\z,-r))\d r \qquad \forall 0 < s < t < \t_-(\z).
\end{equation*}
\end{propo}
\begin{proof} Define, for any $n \geq 1$, $f_n=\varrho_n \diamond
f$, so that, from Propositions \ref{moll2} and \ref{approx},
$\lim_{n \to \infty}\|f_n -f\|  + \|\T_\mathrm{max} f_n
-\T_\mathrm{max} f\| =0.$ In particular, $$\lim_{n \to
\infty}\int_{\O_-}\left|f_n(\x)-f(\x)\right|+\left|\,[\T_\mathrm{max}
f_n](\x)-[\T_\mathrm{max} f ](\x)\right|\d\mu(\x)=0.$$
 Then
Eq. \eqref{10.47} yields
\begin{multline*}
\int_{\Gamma_-}\d\mu_-(\y)\int_0^{\t_+(\y)}\left|f_n({\mathbf{\Phi}}(\y,s))-f({\mathbf{\Phi}}(\y,s))\right|\d
s \\
+ \int_{\Gamma_-}\d\mu_-(\y)\int_0^{\t_+(\y)}\left|[\T_\mathrm{max}
f_n](\mathbf{\Phi}(\y,s))- [\T_\mathrm{max}
f]({\mathbf{\Phi}}(\y,s))\right|\d s \underset{n\to
\infty}{\longrightarrow} 0\end{multline*} since $\T_\mathrm{max} f$
and $\T_\mathrm{max} f_n$ both belong to $L^1(\O,\d\mu)$.
Consequently, for almost every $\y \in \Gamma_-$ (up to a
subsequence, still denoted by $f_n$) we get
\begin{equation*}
\begin{cases}
f_n({\mathbf{\Phi}}(\y,\cdot)) \longrightarrow f({\mathbf{\Phi}}(\y,\cdot))\\
\T_\mathrm{max} f_n ({\mathbf{\Phi}}(\y,\cdot )) \longrightarrow
[\T_\mathrm{max} f ]({\mathbf{\Phi}}(\y,\cdot )) \quad \text{ in }
\quad L^1((0,\t_+(\y))\,,\d s)\end{cases}\end{equation*} as $n \to
\infty$. Let us fix $\y \in \Gamma_-$ for which this holds. Passing
again to a subsequence, we may assume that $f_n({\mathbf{\Phi}}(\y,s
))$ converges (pointwise) to $f({\mathbf{\Phi}}(\y,s ))$ for almost
every $s \in (0,\t_+(\y))$. Let us fix such a $s_0$. Then,
$$f_n({\mathbf{\Phi}}(\y,s_0))-f_n({\mathbf{\Phi}}(\y,s ))=\int_{s_0}^s[\T_\mathrm{max} f_n]({\mathbf{\Phi}}(\y,r ))\d
r \qquad \forall s \in (0,\t_+(\y)).$$ Now, the right-hand-side has
a limit as $n \to \infty$ so that the first term on the left-hand
side also must converge as $n \to \infty.$ Thus, for any $s \in
(0,\t_+(\y))$, the limit
$$\lim_{n \to \infty}f_n({\mathbf{\Phi}}(\y,s))=\widetilde{f}_-({\mathbf{\Phi}}(\y,s ))$$
exists and, for any $0 < s < \t_+(\y)$
\begin{equation*}\label{d/ds}
\widetilde{f}_-({\mathbf{\Phi}}(\y,s
))=\widetilde{f}_-({\mathbf{\Phi}}(\y,s_0 ))-\int_{s_0}^s
[\T_\mathrm{max} f]({\mathbf{\Phi}}(\y,r ))\d r.\end{equation*}  It
is easy to check then that $\widetilde{f}_-(\x)=f(\x)$ for almost
every $\x \in \O_-$. The same arguments lead to the existence of
$\widetilde{f}_+$.
\end{proof}
The above result shows that the mild formulation of Theorem
\ref{representation} is fulfilled for any $\x \in \O_- \cup \O_+$.
It remains to deal with $\O_\infty:=\O_{-\infty}\cap \O_{+\infty}.$
\begin{propo}\label{f+-infinity} Let $f \in \D(\T_\mathrm{max})$. Then, there exists a set $\mathcal{O}
\subset \O_\infty$ with $\mu(\mathcal{O})=0$ and a function
$\widetilde{f}$ defined on $\{\z={\mathbf{\Phi}}(\x,t),\,\x \in
\O_\infty \setminus \mathcal{O},\,t \in \mathbb{R}\}$ such that
$f(\x)=\widetilde{f}(\x)$ $\mu$-almost every $\x \in \O_\infty$ and
$$\widetilde{f}({\mathbf{\Phi}}(\x,s))-\widetilde{f}({\mathbf{\Phi}}(\x,t))=\int_s^t
[\T_\mathrm{max} f ]({\mathbf{\Phi}}(\x, r))\d r, \qquad \forall\,
\x \in \O_\infty \setminus \mathcal{O},\: s < t.$$
\end{propo}
\begin{proof} Since $(\x,t) \mapsto
(\z,t)=(\mathbf{\Phi}(\x,t),t)$ is a measurable and measure
preserving mapping from $\O_\infty \times \mathbb{R}$ onto itself,
Propositions \ref{approx} and \ref{moll2} give
\begin{eqnarray}&&\lim_{n \to
\infty}\int_{\O_\infty}\d\mu(\x)\int_{I_k}\left|f_n(\mathbf{\Phi}(\x,t))-f(\mathbf{\Phi}(\x,t))\right|dt=0\label{jeden}\\
&&\lim_{n \to \infty}
\int_{\O_\infty}\d\mu(\x)\int_{I_k}\left|\T_\mathrm{max}
f_n(\mathbf{\Phi}(\x,t))-\T_\mathrm{max}
f(\mathbf{\Phi}(\x,t))\right|\d t =0,\label{dwa} \end{eqnarray} for
any $I_k=[-k,k]$, $k \in \mathbb{N}$. This shows, in particular,
that there is (a maximal) $\mathcal{E} \subset \O_\infty$ with
$\mu(\mathcal{E})=0$ such that, for almost every $\x \in \O_\infty
\setminus \mathcal{E}$ and any bounded interval $I \subset
\mathbb{R}$:
$$\int_{I} |f(\mathbf{\Phi}(\x,t))|\d t + \int_{I}|[\T_\mathrm{max} f ]
(\mathbf{\Phi}(\x,t))|\d t < \infty$$ and we can argue as in
Proposition \ref{prop312} that
$$f_n(\mathbf{\Phi}(\x,s))-f_n(\x)=-\int_0^s \T_\mathrm{max}
f_n(\mathbf{\Phi}(\x,r))\d r, \qquad \forall s \in \mathbb{R}.$$
Proposition \ref{approx} yields the existence of a subsequence
$(f_{n_p})_p$ and a $\mu$-null set $A_0$ with $\mathcal{E} \subset
A_0 \subset \O_\infty$ such that
$$\lim_{p \to \infty}f_{n_p}(\x)=f(\x), \qquad \forall \x \in
\mathcal{\O_\infty} \setminus \mathcal{A}_0.$$ Now, for any $k \in
\mathbb{N}$,
$$\lim_{p \to
\infty} \int_{\O_\infty}\d\mu(\x)\int_{I_k}\left|\T_\mathrm{max}
f_{n_p}(\mathbf{\Phi}(\x,t))-\T_\mathrm{max}
f(\mathbf{\Phi}(\x,t))\right|\d t =0$$ so that, there is a
subsequence (depending on $k$) and a $\mu$-null set $A_k$ with $A_0
\subset A_k \subset \O_\infty$ such that
$$\lim_{p_{(k)} \to\infty}\int_{I_k}\left|\T_\mathrm{max} f_{n_{p_{(k)}}}(\mathbf{\Phi}(\x,t))-\T_\mathrm{max}
f(\mathbf{\Phi}(\x,t))\right|\d t=0, \qquad \forall \x \in
\O_\infty\setminus A_k.$$ Let $\x \in \O_\infty \setminus A_k$ and
$|s| < k$ be fixed. From
$$f_{n_{p_{(k)}}}(\mathbf{\Phi}(\x,s))-f_{n_{p_{(k)}}}(\x)=-\int_0^s \T_\mathrm{max}
f_{n_{p_{(k)}}}(\mathbf{\Phi}(\x,r))\d r$$ we deduce that the limit
$\lim_{p_{(k)} \to \infty}f_{n_{p_{(k)}}}(\mathbf{\Phi}(\x,s))$
exists and is equal to  $$\lim_{p_{(k)} \to
\infty}f_{n_{p_{(k)}}}(\mathbf{\Phi}(\x,s))=f(\x)-\int_0^s
\T_\mathrm{max}
 f(\mathbf{\Phi}(\x,r))\d r.$$
 We define then $\widetilde{f}$ by
 $$\widetilde{f}(\mathbf{\Phi}(\x,s))=\lim_{p_{(k)} \to
\infty}f_{n_{p_{(k)}}}(\mathbf{\Phi}(\x,s)), \qquad \x \in \O_\infty
\setminus A_k,\:|s| < k$$ and defining $\mathcal{O}=\bigcup_{k \geq
1} A_k$, we get the result.
\end{proof}

Before the proof of Theorem \ref{representation}, we have to
establish existence of the trace on $\Gamma_-$.
\begin{propo}\label{trace} Let $f$ satisfies condition (1) of Theorem \ref{representation}.
Then  $$\lim_{t \to 0+}f^\sharp({\mathbf{\Phi}}(\y,t))$$ exists  for
almost every $\y \in \Gamma_-$. Similarly, $\lim_{t \to
0+}f^\sharp({\mathbf{\Phi}}(\y,-t))$ exists for almost every $\y \in
\Gamma_+$.
\end{propo}
\begin{proof} First we note that there is $\widetilde{\Omega}_-\subset \Omega_-$ with
 $\mu(\Omega_-\setminus\widetilde{\Omega}_-)=0$ such that (\ref{integralTM-}) is valid any $\x \in \widetilde{\Omega}_-$. Let
$\widetilde{\Gamma}_-=\{\y\in \Gamma_-;\; \y =
\mathbf{\Phi}(\x,\tau_-(\x)), \x \in \widetilde{\Omega}_-\}$. It is
easy to see that $\mu_-(\Gamma_-\setminus \widetilde{\Gamma}_-) =0$.
Indeed, otherwise, by (\ref{10.47}), there would be a subset of
$\Omega_-$ of positive $\mu$-measure, not intersecting
$\widetilde{\Omega}_-$, which would contradict (\ref{integralTM-}).
Consequently,  any $\x\in \widetilde{\Omega}_-$ can be written as
$\x = \Phi(\y, \tau_-(\y))$, $\y \in \widetilde{\Gamma}_-$ and
(\ref{integralTM-}) can be recast as
\begin{equation}
f^\sharp({\mathbf{\Phi}}(\y,t))-f^{\sharp}({\mathbf{\Phi}}(\y,t_0))=\int_{t}^{t_0}g(\mathbf{\Phi}(\y,s)\d
s. \label{nowe}
\end{equation}
 for almost any $\y \in \Gamma_-$, where $0<t\leq t_0<
\tau_+(\y)$. Using again (\ref{10.47}), $s\mapsto
g({\mathbf{\Phi}}(\y,s)$ is integrable on $(0,\tau_+(\y)$ for almost
any $\y\in \Gamma_-$. Consequently, for almost every $\y\in
\Gamma_-$ we can pass to the limit in (\ref{nowe}) with $t\to 0$; it
is easy to check that this limit does not depend on $t_0$. The
existence of $\lim_{t\to 0+}f^\sharp({\mathbf{\Phi}}(\y,-t))$ for a.
e. $\y \in \Gamma_+$ follows by the same argument.
\end{proof}

The above proposition allows to define the trace operators.
\begin{defi}
For any $f \in \D(\T_{\mathrm{max}})$, define the \textit{traces}
$\B^{\pm}f$ by
\begin{equation*}
\B^+f(\y):=\lim_{t \to 0+}f^\sharp({\mathbf{\Phi}}(\y,-t)) \qquad
\text{ and } \qquad \B^-f(\y):=\lim_{t \to
0+}f^\sharp({\mathbf{\Phi}}(\y,t))
\end{equation*} for any $\y \in \Gamma_\pm$ for which the limits
exist, where   $f^\sharp$ is a suitable representative of $f$.
\end{defi}

\begin{preuve2}
To prove that $(2) \implies (1)$,  given $f \in
\D(\T_\mathrm{max})$, set
\begin{equation*}
f^\natural(\x)=\begin{cases} \widetilde{f}_-(\x)  \quad &\text{ if }
\x
\in \O_-,\\
 \widetilde{f}_+(\x)  \quad &\text{ if } \x
\in \O_+ \cap \O_{-\infty},\\
 \widetilde{f}(\x) \quad &\text{ if } \x
\in \O_{-\infty} \cap \O_{+\infty},\\
\end{cases}
\end{equation*}
where $\widetilde{f}_{\pm}$ are given by Proposition \ref{f+-} while
$\widetilde{f}$ is provided by Prop. \ref{f+-infinity}. Then, it is
clear that for any $\x \in
  \O$ and any $-\t_-(\x) < t_1 \leq t_2 <
  \t_+(\x)$
  $$f^\sharp({\mathbf{\Phi}}(\x,t_1))-f^\sharp({\mathbf{\Phi}}(\x,t_2))=\int_{t_1}^{t_2} [\T_\mathrm{max} f ]({\mathbf{\Phi}}(\x,s))\d
  s$$
and \eqref{integralTM-} is proven.

Let us now prove that $(1) \implies (2)$. Let us fix $\psi \in
\mathfrak{Y}$, one has
\begin{equation*}\begin{split}
\int_{\O_-} f(\x) \dfrac{\d}{\d s}
\psi(\mathbf{\Phi}(\x,s))\big|_{s=0}\d\mu(\x)&=\int_{\Gamma_-}\d\mu_-(\y)\int_0^{\t_+(\y)}
 f(\mathbf{\Phi}(\y,t))  \dfrac{\d}{\d t}\psi(\mathbf{\Phi}(\y,t))\d t\\
&=\int_{\Gamma_-}\d\mu_-(\y)\int_0^{\t_+(\y)}
f^\sharp(\mathbf{\Phi}(\y,t)) \dfrac{\d}{\d
t}\psi(\mathbf{\Phi}(\y,t))\d t.
\end{split}
\end{equation*}
Notice that since both $\int_{\O_-}f(\x) \dfrac{\d}{\d s}
\psi(\mathbf{\Phi}(\x,s))\big|_{s=0}\d\mu(\x)$ and $\int_{\O_-}
\psi(\x)g(\x)\d\mu(\x)$ exist, Proposition \ref{prointegra} and
Fubini's Theorem, the integrals $$\int_0^{\t_+(\y)}
f^\sharp(\mathbf{\Phi}(\y,t)) \dfrac{\d}{\d
t}\psi(\mathbf{\Phi}(\y,t))\d t \quad \text{ and } \quad
\int_0^{\t_+(\y)} g(\mathbf{\Phi}(\y,t)) \psi(\mathbf{\Phi}(\y,t))\d
t$$ are well-defined for $\mu_-$-almost every $\y \in \Gamma_-$. Let
us prove that these two integrals coincide for almost-every $\y \in
\Gamma_-$. According to Lemma \ref{supp}, for almost every $\y \in
\Gamma_-$, there is a sequence $(t_n)_n$ (depending on $\y$) such
that $\psi(\mathbf{\Phi}(\y,t_n))=0$ and $t_n \to \t_+(\y)$. Thus,
$$\int_0^{\t_+(\y)}
f^\sharp(\mathbf{\Phi}(\y,t)) \dfrac{\d}{\d
t}\psi(\mathbf{\Phi}(\y,t))\d t=\lim_{n \to \infty}\int_0^{t_n}
f^\sharp(\mathbf{\Phi}(\y,t)) \dfrac{\d}{\d
t}\psi(\mathbf{\Phi}(\y,t))\d t$$ and
$$\int_0^{\t_+(\y)} g(\mathbf{\Phi}(\y,t)) \psi(\mathbf{\Phi}(\y,t))\d
t=\lim_{n \to
\infty}\int_0^{t_n}\psi(\mathbf{\Phi}(\y,t))g(\mathbf{\Phi}(\y,t))\d
t.$$ Further, for almost every $\y \in \Gamma_-$, according to
\eqref{integralTM-},
$$f^\sharp(\mathbf{\Phi}(\y,t))=\B^-f(\y)-\int_{0}^t
g(\mathbf{\Phi}(\y,r))\d r, \qquad \forall t \in (0,\t_+(\y)).$$
Integration by parts, using the fact that
$\psi(\mathbf{\Phi}(\y,0))=\psi(\mathbf{\Phi}(\y,t_n))=0$ for any
$n$,  leads to
$$\int_0^{t_n}
f^\sharp(\mathbf{\Phi}(\y,t)) \dfrac{\d}{\d
t}\psi(\mathbf{\Phi}(\y,t))\d
t=\int_0^{t_n}g(\mathbf{\Phi}(\y,t))\psi(\mathbf{\Phi}(\y,t))\d t.$$
Consequently, for $\mu_-$ almost every $\y \in \Gamma_-$:
\begin{equation}\label{tfini}\int_{0}^{\t_+(\y)}
f^\sharp(\mathbf{\Phi}(\y,t)) \dfrac{\d}{\d
t}\psi(\mathbf{\Phi}(\y,t))\d
t=\int_{0}^{\t_+(\y)}\psi(\mathbf{\Phi}(\y,s))g(\mathbf{\Phi}(\y,t))
\d t.\end{equation} Finally, we get
\begin{equation}\label{og-}\begin{split}
\int_{\O_-} f(\x) \dfrac{\d}{\d s}
\psi(\mathbf{\Phi}(\x,s))\big|_{s=0}\d\mu(\x)&=\int_{\Gamma_-}\d\mu_-(\y)\int_0^{\t_+(\y)}\psi(\mathbf{\Phi}(\y,t))g(\mathbf{\Phi}(\y,t))\d
t \\
&=\int_{\O_-} g( \x) \psi(\x)\d\mu(\x).
\end{split}
\end{equation}  Using
now parametrization over $\Gamma_+$, we prove in the same way  that
\begin{equation}\label{og+}\int_{\O_+ \cap \O_{-\infty}} f(\x) \dfrac{\d}{\d s}
\psi(\mathbf{\Phi}(\x,s))\big|_{s=0}\d\mu(\x)=\int_{\O_+ \cap
\O_{-\infty}}g(\x) \psi(\x)\d\mu(\x).\end{equation} It remains now
to evaluate $A:=\int_{\O_{+\infty} \cap \O_{-\infty}} f(\x)
\frac{\d}{\d s} \psi(\mathbf{\Phi}(\x,s))\big|_{s=0}\d\mu(\x).$
According to Assumption \ref{ass:h2}  $$A =\int_{\O_{+\infty} \cap
\O_{-\infty}}
 f^\sharp(\mathbf{\Phi}(\x,t))  \dfrac{\d}{\d t}
\psi(\mathbf{\Phi}(\x,t))\d\mu(\x), \qquad \forall t \in
\mathbb{R}.$$ Let us integrate the above identity over $(0,1)$, so
that
$$  A =\int_{\O_{-\infty} \cap \O_{+\infty}}\d\mu(\x)\int_0^1 f^\sharp(\mathbf{\Phi}(\x,t))  \dfrac{\d}{\d t}
\psi(\mathbf{\Phi}(\x,t)) \d t.$$ Let us fix $\x \in \O_{-\infty}
\cap \O_{+\infty}$. For any $t \in (0,1)$, one  has
 $f^\sharp(\mathbf{\Phi}(\x,t))=f^\sharp(\x)-\int_0^t g(\mathbf{\Phi}(\x,s))\d s$
and integration by parts yields
\begin{multline*}
\int_0^1 f^\sharp(\mathbf{\Phi}(\x,t))\dfrac{\d }{\d
t}\psi(\mathbf{\Phi}(\x,t))\d
t=\int_0^1 \psi(\mathbf{\Phi}(\x,t))g(\mathbf{\Phi}(\x,t))\d t -\psi(\x)f^\sharp(\x)\\
+\psi(\mathbf{\Phi}(\x,1))\bigg(f^\sharp(\x)
 -\int_0^1
g(\mathbf{\Phi}(\x,s))\d s\bigg)\\
=\int_0^1 \psi(\mathbf{\Phi}(\x,t))g(\mathbf{\Phi}(\x,t))\d t
+\psi(\mathbf{\Phi}(\x,1))f^\sharp(\mathbf{\Phi}(\x,1))-\psi(\x)f^\sharp(\x)
\end{multline*}
where we used again \eqref{integralTM-}. Integrating over
$\O_{-\infty} \cap \O_{+\infty}$ we see from Liouville's Theorem
(Assumption \ref{ass:h2}) that
$$\int_{\O_{-\infty} \cap \O_{+\infty}}
 \psi(\mathbf{\Phi}(\x,1))f^\sharp(\mathbf{\Phi}(\x,1))\d\mu(\x)=\int_{\O_{-\infty} \cap
 \O_{+\infty}}\psi(\x)f^\sharp(\x)\d\mu(\x),$$ i.e.
$$A=\int_{\O_{-\infty} \cap \O_{+\infty}}\d\mu(\x)\int_0^1\psi(\mathbf{\Phi}(\x,t))g(\mathbf{\Phi}(\x,t))\d
t$$ which, thanks to Liouville's Theorem, is nothing but
\begin{equation}\label{og+infty}
\int_{\O_{+\infty} \cap \O_{-\infty}} f(\x) \frac{\d}{\d s}
\psi(\mathbf{\Phi}(\x,s))\big|_{s=0}\d\mu(\x)=\int_{\O_{-\infty}
\cap \O_{+\infty}} g(\x)\,\psi(\x)\d\mu(\x).\end{equation} Combining
\eqref{og-}, \eqref{og+} and \eqref{og+infty}, we obtain
$$\int_\O f(\x)\frac{\d}{\d s}
\psi(\mathbf{\Phi}(\x,s))\big|_{s=0}\d\mu(\x)=\int_\O g(\x)
\psi(\x)\d\mu(\x), \qquad \forall \psi \in \mathfrak{Y}$$ which
exactly means that $f \in \D(\T_\mathrm{max})$ with
$g=\T_\mathrm{max}$ and the proof is complete.\end{preuve2}

\begin{cor}\label{cortrace} Traces $B^{\pm}f$ on $\Gamma_\pm$ can be defined for any
$f \in \D(\T_{\mathrm{max}})$. For $\mu_-$- almost any $\y \in
\Gamma_-$ we have
$$\B^-f(\y)=f^\sharp({\mathbf{\Phi}}(\y,t))+\int_0^t[\T_{\mathrm{max}}f]({\mathbf{\Phi}}(y, s))\d
s, \qquad \forall t \in (0,\tau_+(y)),$$ where $f^\sharp$ is a
suitable representative of $f$. An analogous formula holds for
$\B^+f$.\end{cor}

Lemma \ref{defmu+-} provides the existence of Borel measures
$\d\mu_{\pm}$ on $\Gamma_{\pm}$, which allow us to define the natural trace
spaces associated to Problem \eqref{1}, namely,
$$L^1_{\pm}:=L^1(\Gamma_\pm,\d\mu_\pm).$$
However, the traces $B^\pm f$, $f \in \D(\T_{\mathrm{max}})$, not
necessarily belong to $L^1_{\pm}.$

\section{Well-posedness for initial and boundary- value problems}

\subsection{Absorption semigroup}\label{sec:exis}
From now on, we will denote $X=L^1(\O,\d\mu)$ endowed with its
natural norm $\|\cdot\|_X$. Let $\T_0$ be the free streaming
operator with \textit{no re--entry boundary conditions}:
$$\T_0\psi=\T_{\mathrm{max}}\psi, \qquad \text{ for any }
\psi \in \D(\T_0),$$ where the domain $\D(\T_0)$ is defined by
$$\D(\T_0)=\{\psi \in \D(\T_{\mathrm{max}})\,;\,\B^-\psi=0\}.$$
We state the following generation result:
\begin{theo}\label{uot} The operator $(\T_0,\D(\T_0))$ is the generator of a nonnegative $C_0$-semigroup
of contractions $\uot$ in $L^1(\O,\d\mu)$ given by
$$U_0(t)f(\x)=f({\mathbf{\Phi}}(\x,-t))\chi_{\{t <
\tau_-(\x)\}}(\x), \qquad (\x \in \O,\:f \in X),$$ where $\chi_A$
denotes the characteristic function of a set $A$.
\end{theo}
\begin{proof} The proof is divided into three steps:

$\bullet$ {\it Step 1.} Let us first check that the family of
operators $\uot$ is a nonnegative contractive $C_0$-semigroup in
$X$. Thanks to Proposition \ref{Phiprop}, we can prove that, for any
$f \in X$ and any $t \geq 0$, the mapping $U_0(t)f\::\O \to
\mathbb{R}$ is measurable and the semigroup properties $U_0(0)f=f$
and $U_0(t)U_0(s)f=U_0(t+s)f$ $(t,s \geq 0)$ hold. Let us now show
that $\|U_0(t)f\|_X \leq \|f\|_X$. We have
$$\|U_0(t)f\|_X=\int_{\O_+}|U_0(t)f|\d\mu+\int_{\O_-\cap
\O_{+\infty}}|U_0(t)f|\d\mu+\int_{\O_{-\infty}\cap\O_{+\infty}}|U_0(t)f|\d\mu.$$
Propositions \ref{prointegra} and \ref{Phiprop} yield
\begin{equation*}\begin{split}
\int_{\O_+}|U_0(t)f|\d\mu
&=\int_{\Gamma_+}\d\mu_+(\y)\int_0^{\tau_-(\y)}
|U_0(t)f({\mathbf{\Phi}}(\y,-s))|\d s\\
&= \int_{\Gamma_+}\d\mu_+(\y)\int_0^{\max(0,\tau_-(\y)-t)}
|f({\mathbf{\Phi}}(\y,-s-t))|\d s\\
&\leq \int_{\Gamma_+}\d\mu_+(\y)\int_t^{\max(t,\tau_-(\y))}
|f({\mathbf{\Phi}}(\y,-r))|\d r \leq \int_{\O_+}|f|\d\mu.
\end{split}\end{equation*}
In the same way we obtain
\begin{equation*}
\int_{\O_- \cap \O_{+\infty}}|U_0(t)f| \d\mu
=\int_{\Gamma_{-\infty}}\d\mu_-(\y)\int_0^{\infty}
|U_0(t)f({\mathbf{\Phi}}(\y, s))|\d s= \int_{\O_- \cap
\O_{+\infty}}|f|\d\mu,
\end{equation*}
and
\begin{equation*}
\int_{\O_{-\infty}\cap \O_{+\infty}}|U_0(t)f|\d\mu =
\int_{\O_{-\infty}\cap \O_{+\infty}}|f|\d\mu.\end{equation*} This
proves contractivity of $U_0(t)$. Let us now show that $U_0(t)f$ is
continuous, i.e.
$$\lim_{t \to 0}\|U_0(t)f-f\|_X=0.$$
It is enough to show that this property holds for any $f \in
\Con_0(\O)$. In this case, $\lim_{t\to 0}U_0(t)f(\x)=f(\x)$ for any
$\x \in \O.$ Moreover, $\sup_{\x \in \O}|U_0(t)f(\x)| \leq \sup_{\x
\in \O}|f(\x)|$ and the support of $U_0(t)f$ is bounded, so that the
Lebesgue dominated convergence theorem leads to the result. This
proves that $\uot$ is a $C_0$-semigroup of contractions in $X$. Let
$\A_0$ denote its generator.

$\bullet$  \textit{Step 2.} To show that $\D(\A_0) \subset
\D(\T_0)$,   fix $f \in \D(\A_0)$, $\lambda > 0$ and $g=(\lambda-
\A_0)f.$ Then,
\begin{equation}f(\x)=\int_0^{\t_-(\x)}\exp(-\lambda
t)\,g({\mathbf{\Phi}}(\x,-t))\d t, \qquad (\x \in \O). \label{ho}
\end{equation}
  To prove
that $f \in \D(\T_\mathrm{max})$ with $\T_\mathrm{max}f=\A_0f$, it
suffices to prove that
$$\int_{\O
}(\lambda f(\x)-g(\x))\psi(\x)\d\mu(\x)=\int_{\O }f(\x)\dfrac{\d}{\d
  s}\psi({\mathbf{\Phi}}(\x,s))\bigg|_{s=0} \d\mu(\x),\qquad \qquad \forall \psi \in
\mathfrak{Y}.$$ Let us fix $\psi \in \mathfrak{Y}$, set
$\varphi(\x):=\frac{\d}{\d
  s}\psi({\mathbf{\Phi}}(\x,s))\big |_{s=0}$ and write
\begin{multline*}
\int_{\O}  f(\x)\varphi(\x)\d\mu(\x)=\int_{\O_+}
f(\x)\varphi(\x)\d\mu(\x)+\int_{\O_{+\infty} \cap \O_-}
f(\x)\varphi(\x)\d\mu(\x)\\+\int_{\O_{+\infty} \cap \O_{-\infty} }
f(\x)\varphi(\x)\d\mu(\x)=I_1+I_2+I_3.
\end{multline*}
We first deal with $I_1$. For any $\y \in \Gamma_+$ and $t \in
(0,\t_-(\y))$ we have $\varphi({\mathbf{\Phi}}(\y,-t))=-\frac{\d}{\d
t}\psi({\mathbf{\Phi}}(\y,-t))$ and
$f({\mathbf{\Phi}}(\y,-t))=\int_t^{\tau_-(\y)}\exp(-\lambda
(s-t))g({\mathbf{\Phi}}(\y,-s))\d s.$  Then, by Proposition
\ref{prointegra},
\begin{equation*}\begin{split}
I_1&=-\int_{\Gamma_+}\d\mu_+(\y)\int_0^{\t_-(\y)}\frac{\d}{\d
t}\psi({\mathbf{\Phi}}(\y,-t))\d t \int_t^{\t_-(\y)}\exp(-\lambda
(s-t))g({\mathbf{\Phi}}(\y,-s))\d s\\
&=-\int_{\Gamma_+}\d\mu_+(\y)\int_0^{\tau_-(\y)}g({\mathbf{\Phi}}(\y,-s))\d
s\int_0^s \exp(-\lambda (s-t))\frac{\d}{\d
t}\left(\psi({\mathbf{\Phi}}(\y,-t ))\right)\d t\\
&=\int_{\Gamma_+}\d\mu_+(\y)\int_0^{\tau_-(\y)}g({\mathbf{\Phi}}(\y,-s )) \times \\
& \phantom{\int_{\Gamma_+}\d\mu_+(\y)}\times\left\{ \lambda \int_0^s
\exp(-\lambda (s-t))\psi({\mathbf{\Phi}}(\y,-t))\d t -
\psi({\mathbf{\Phi}}(\y,-s))\right\}\d s\end{split}
\end{equation*}
where we used the fact that $\psi(\mathbf{\Phi}(\y,0))=0$ for
 any $\y \in \Gamma_+$ since $\psi$ is compactly supported.
Thus
\begin{equation*}\begin{split}
I_1&=\lambda\int_{\Gamma_+}\d\mu_+(\y)\int_0^{\t_-(\y)}\psi(\mathbf{\Phi}(\y,-t))\d
t \int_t^{\t_-(\y)}\exp(-\lambda(s-t))g(\mathbf{\Phi}(\y,-s))\d s
\\
&\phantom{++++++++++++}
-\int_{\Gamma_-}\d\mu_+(\y)\int_0^{\t_-(\y)}g(\mathbf{\Phi}(\y,-s))\psi(\mathbf{\Phi}(\y,-s))\d
s\\
&=
\int_{\Gamma_+}\d\mu_+(\y)\int_0^{\t_-(\y)}\psi(\mathbf{\Phi}(\y,-t))\big(\lambda
f(\mathbf{\Phi}(\y,-t))-g(\mathbf{\Phi}(\y,-t))\big)\d t.
\end{split}
\end{equation*}
 Using again
Proposition \ref{prointegra}, we obtain
\begin{equation}\label{I1}I_1=\int_{\O_+}\left(\lambda f(\x)-g(\x)\right)\psi(\x)\d\mu(\x).\end{equation}
Arguing in a similar way, we prove that
\begin{equation}\label{I2}
I_2=-\int_{\O_- \cap \O_{+\infty}}\left(\lambda
f(\x)-g(\x)\right)\psi(\x)\d\mu(\x).\end{equation} Finally, since
$$f(\x)=\int_0^{\infty}\exp(-\lambda t)g\left({\mathbf{\Phi}}(\x,-t)\right)\d t
\qquad \text{ for any } \quad\x \in \O_{-\infty} \cap
\O_{+\infty},$$ one has
\begin{equation*}\begin{split}
I_3&=\int_{\O_{-\infty} \cap
\O_{+\infty}}\varphi(\x)\d\mu(\x)\int_0^\infty
\exp(-\lambda t) g({\mathbf{\Phi}}(\x,-t))\d t\\
&=\int_0^\infty \exp(-\lambda t)\d t \int_{\O_{-\infty} \cap
\O_{+\infty}}\varphi(\x)g({\mathbf{\Phi}}(\x,-t))\d\mu(\x).\end{split}\end{equation*}
Now, Assumption \ref{ass:h2}  asserts that
$$\int_{\O_{-\infty} \cap
\O_{+\infty}}\varphi(\x)g({\mathbf{\Phi}}(\x,-t))\d\mu(\x)=\int_{\O_{-\infty}
\cap \O_{+\infty}}g(\x)\varphi({\mathbf{\Phi}}(\x, t))\d\mu(\x),
\qquad \forall t \geq 0,$$ and, since
$\varphi({\mathbf{\Phi}}(\x,t))=\frac{\d }{\d t}
\psi({\mathbf{\Phi}}(\x,t))$,  finally
\begin{equation*}\begin{split}
I_3&=\int_{\O_{-\infty} \cap
\O_{+\infty}}g(\x)\d\mu(\x)\int_0^{\infty}\exp(-\lambda
t)\dfrac{\d}{\d
t}\left(\psi({\mathbf{\Phi}}(\x ,t))\right)\d t\\
&=-\int_{\O_{-\infty} \cap \O_{+\infty}}g(\x)\psi(\x)\d\mu(\x) +
\lambda \int_{\O_{-\infty} \cap
\O_{+\infty}}g(\x)\d\mu(\x)\int_0^\infty \exp(-\lambda
t)\psi({\mathbf{\Phi}}(\x, t))\d t.\end{split}\end{equation*} Using
again Assumption \ref{ass:h2}, this finally gives
\begin{equation}\label{I3}
I_3=-\int_{\O_{-\infty} \cap \O_{+\infty}}\left(g(\x)-\lambda
f(\x)\right)\psi(\x)\d\mu(\x).\end{equation} Combining
\eqref{I1}--\eqref{I3} leads to
$$\int_{\O }f(\x)\dfrac{\d}{\d
  s}\psi({\mathbf{\Phi}}(\x,s))\bigg|_{s=0} \d\mu(\x)=-\int_{\O}\left(g(\x)-\lambda
f(\x)\right)\psi(\x)\d\mu(\x)$$ which proves that $f \in
\D(\T_{\mathrm{max}})$ and $(\lambda- \T_{\mathrm{max}})f=g.$ Next,
for $\y \in \Gamma_-$ and $0 < t < \tau_+(\y)$ we write
$t=\tau_-\left({\mathbf{\Phi}}(\y,t)\right)$ and, by Proposition
\ref{Phiprop} and (\ref{ho}), we obtain
\begin{equation}\label{fy+}
f({\mathbf{\Phi}}(\y,t)) = \int_0^t\exp(-\lambda
(t-s))\,g({\mathbf{\Phi}}(\y,s))\d s.
\end{equation} Consequently, $\lim_{t \to
0^+}f({\mathbf{\Phi}}(\y,t))=0$  a.e. $\y \in \Gamma_-,$ i.e.
$\B^-f=0$ so that $f \in \D(\T_0)$ and $\A_0f=\T_0f=\lambda f-g$.

 $\bullet$ \textit{Step 3.} Now let us show the converse
inclusion $\D(\T_0) \subset \D(\A_0)$. Let $f \in \D(\T_0)$.
Changing possibly $f$ on a set of zero measure, we may write
$f=f^\sharp$, where $f^\sharp$ is the representative of $f$ given by
Theorem \ref{representation}. Then, for any $\x \in \O$ and any $0
\leq t<\t_-(\x)$
$$f({\mathbf{\Phi}}(\x,-t))-f(\x)=\int_0^t
[\T_\mathrm{max}f]({\mathbf{\Phi}}(\x,-r))\d r$$ which, according to
the explicit expression of $U_0(t)$,  means that
\begin{equation}\label{uot0}
U_0(t)f(\x)-f(\x)=\int_0^t\,U_0(r)\T_{\mathrm{max}}f(\x)\d
r\end{equation}  for any  $\x \in \O$ and $t < \tau_-(\x).$  Letting
$t$ converge towards $\t_-(\x)$ we obtain
$$f(\x)=-\int_0^{\t_-(\x)}[\T_\mathrm{max}f]({\mathbf{\Phi}}(\x,-r))\d r.$$
In particular,  Eq. \eqref{uot0} holds true for any $\x \in \O$ and
any $t \geq \tau_-(\x).$ Arguing exactly as in \cite[p. 38]{voigt},
the pointwise identity \eqref{uot0} represents the $X$--integral,
i.e, $U_0(t)f-f=\int_0^t \,U_0(r)\T_{\mathrm{max}}f \d r$ in
$L^1(\O,\d\mu)$. Consequently, $f \in \D(\A_0)$ with
$\A_0f=\T_{\mathrm{max}}f.$\end{proof}

\subsection{\textbf{Green's formula}}The above result allows us to treat more general boundary-value
problem:
\begin{theo}\label{Theo10.43} Let $u \in \lm$ and $g \in X$ be
given. Then the function
\begin{equation*}
f(\x)=\int_0^{\t_-(\x)}\exp(-\lambda t)\,g({\mathbf{\Phi}}(\x,-t))\d
t +\chi_{\{\t_-(\x) < \infty\}} \exp(-\lambda \t_-(\x))
u({\mathbf{\Phi}}(\x,-\t_-(\x)))\end{equation*} is a \textbf{
unique} solution $f \in \D(\T_{\mathrm{max}})$ of the boundary value
problem:
\begin{equation}\label{BVP1}
\begin{cases}
(\lambda- \T_{\mathrm{max}})f=g,\\
\B^-f=u,
\end{cases}
\end{equation}
where $\lambda > 0.$ Moreover, $\B^+f \in \lp$ and
\begin{equation}\label{10.1156} \|\B^+f\|_{\lp} + \lambda \|f\|_X \leq
\|u\|_{\lm} + \|g\|_X,\end{equation} with equality sign if $g \geq
0$ and $u \geq 0$.
\end{theo}
\begin{proof} Let us write $f=f_1+f_2$ with $f_1(\x)=\int_0^{\t_-(\x)}\exp(-\lambda t)\,g({\mathbf{\Phi}}(\x,-t))\d
t,$ and
$$f_2(\x)=\chi_{\{\t_-(\x) < \infty\}} \exp(-\lambda \t_-(\x)) u\big({\mathbf{\Phi}}\left(\x,-\t_-(\x)\right)\big), \qquad \x \in
\O.$$ According to Theorem \ref{uot}, $f_1=(\lambda- \T_0)^{-1}g,$
i.e. $f_1 \in \D(\T_{\mathrm{max}})$ with $(\lambda-
\T_{\mathrm{max}})f_1=g$ and $\B^-f_1=0$. Therefore, to prove that
$f$ is a solution of \eqref{BVP1} it suffices to check that $f_2 \in
\D(\T_{\mathrm{max}})$, $(\lambda- \T_{\mathrm{max}})f_2=0$ and
$\B^-f_2=u.$ It is easy to see that $f_2 \in L^1(\O,\d\mu)$ (see
also \eqref{f-2uu}). To prove that $f_2 \in \D(\T_\mathrm{max})$ one
argues as in the proof of Theorem \ref{uot}. Precisely, let $\psi
\in \mathfrak{Y}$,
 noticing that $f_2$ vanishes outside $\O_-$, one has thanks to
 \eqref{f2-u}
\begin{equation*}\begin{split}
\int_{\O}f_2(\x)\dfrac{\d }{\d
s}\psi({\mathbf{\Phi}}(\x,s))\big|_{s=0}\d\mu(\x)
&=\int_{\Gamma_-}\d\mu_-(\y)\int_0^{\t_+(\y)}f_2({\mathbf{\Phi}}(\y,t))\frac{\d
}{\d t} \psi({\mathbf{\Phi}}(\y,t))\d t\\
&=\int_{\Gamma_-} u(\y)\d\mu_-(\y)\int_0^{\t_+(\y)} \exp(-\lambda t)
\frac{\d }{\d t} \psi({\mathbf{\Phi}}(\y,t))\d t.
\end{split}
\end{equation*}
For almost every $\y \in \Gamma_-$, we compute the integral over
$(0,\t_+(\y))$ by parts, which yields $f_2 \in \D(\T_\mathrm{max} )$
with $\T_\mathrm{max}f_2=\lambda f_2.$ Also,
\begin{equation}\label{f2-u}f_2({\mathbf{\Phi}}(\y,t))=\exp(-\lambda t)u(\y), \qquad  \y \in
\Gamma_-,\;\;0<t<\tau_+(\y)\end{equation} from which we see that
$\B^- f_2=u.$

Consequently, $f$ is a solution to \eqref{BVP1}. To prove that the
solution is unique, it is sufficient to prove that the only solution
$h \in \D(\T_{\mathrm{max}})$ to
 $(\lambda- \T_{\mathrm{max}})h=0, \B^-h=0, $
is $h=0$. This follows from the fact that such a solution $h$
actually belongs to $\D(\T_0)$ if $\lambda \in \varrho(\T_0)$.
Finally, it remains to prove \eqref{10.1156}. For simplicity, we
denote the representative of $f_i$, $i=1,2$, defined in Proposition
\ref{trace}, with the same letter. Using \eqref{f2-u} and the fact
that $f_2$ vanishes on $\O_{-\infty}$,  from \eqref{10.47} we get
\begin{equation}\label{f-2uu}\begin{split}
\lambda \int_{\O}|f_2|\d\mu&=\lambda
\int_{\O_-}|f_2|\d\mu=\lambda \int_{\Gamma_-}\d\mu_-(\y)\int_0^{\tau_+(\y)}e^{-\lambda t}|u(\y)|\d t\\
&=\int_{\Gamma_-}|u(\y)|\left(1-e^{-\lambda
\tau_+(\y)}\right)\d\mu_-(\y).\end{split}\end{equation} Define $h$ :
$\y \in \Gamma_- \longmapsto h(\y)=|u(\y)|e^{-\lambda \tau_+(\y)}.$
It is clear that $h$ vanishes on $\Gamma_{-\infty}$ and $h(\y) \leq
|u(\y)|$ for a.e. $\y \in \Gamma_-$. In particular, $h \in \lm$ and,
according to \eqref{10.51},
\begin{equation*}\begin{split}
\int_{\Gamma_-}&h(\y)\d\mu_-(\y) =\int_{\Gamma_- \setminus
\Gamma_{-\infty}}h(\y)\d\mu_-(\y)=\int_{\Gamma_+ \setminus
\Gamma_{+\infty}}h({\mathbf{\Phi}}(z,-\tau_-(z)))\d\mu_+(z)\\
&=\int_{\Gamma_+\setminus \Gamma_{+\infty}}e^{-\lambda
\tau_-(z)}|u({\mathbf{\Phi}}(z,-\tau_-(z)))|\d\mu_+(z)=\int_{\Gamma_+}|\B^+f_2(z)|\d\mu_+(z)=\|\B^+f_2\|_{\lp}.
\end{split}\end{equation*}
Combining this with \eqref{f-2uu} leads to
\begin{equation}\label{B+f-2}\lambda \|f_2\|_X+\|\B^+f_2\|_{\lp}=\|u\|_{\lm}.\end{equation} Now, let us show that $\B^+f_1 \in
\lp$ and $\|\B^+f_1\|_{\lp} + \lambda \|f_1\|_X \leq \|g\|_X.$ For
any $\y \in \Gamma_+$ and $0 < t < \tau_-(\y),$ we see, as above,
that $f_1({\mathbf{\Phi}}(\y,-t))=\int_t^{\tau_-(\y)}\exp(-\lambda
(s-t))g({\mathbf{\Phi}}(\y,-s))\d s.$ This shows that
$$\B^+f_1(\y)=\lim_{t \to
0^+}f_1({\mathbf{\Phi}}(\y,-t))=\int_0^{\tau_-(\y)}\exp(-\lambda
s))g({\mathbf{\Phi}}(\y,-s))\d s.$$ According to Proposition
\ref{prointegra},
$$\int_{\Gamma_+}\d\mu_+(\y)\int_0^{\tau_-(\y)}\left|g({\mathbf{\Phi}}(\y,-s))\right|\d
s=\int_{\O_+}|g|\,\d\mu$$ which, since $\exp(-\lambda
(s-t))|g({\mathbf{\Phi}}(\y,-s))| \leq |g({\mathbf{\Phi}}(\y,-s))|,$
implies $\B^+f_1 \in \lp.$ Let us now assume $g \geq 0$. Then $f_1
\geq 0$ and hence
$$\lambda \|f_1\|=\lambda \int_{\O}f_1\,\d\mu=\lambda \int_{\O_+}f_1\,\d\mu+\lambda \int_{\O_-\cap \O_{+\infty}}f_1\,\d\mu
+\lambda \int_{\O_{-\infty} \cap \O_{+\infty}}f_1\,\d\mu.$$ Using
similar arguments to those used in the study of $f_2$, we have
$$\lambda \int_{\O_+}f_1\,\d\mu=\int_{\Gamma_+}\d\mu_+(\y)\int_0^{\tau_-(\y)}g({\mathbf{\Phi}}(\y,-t))\left(1-\exp(-\lambda t)\right)\d t,$$
 which,  by Proposition \ref{prointegra}, implies
 $\lambda
 \int_{\O_+}f_1\,\d\mu=\int_{\O_+}g\,\d\mu-\int_{\Gamma_+}\B^+f_1\,\d\mu_+.$
Similar argument shows that $\lambda \int_{\O_- \cap
\O_{+\infty}}f_1\,\d\mu= \int_{\O_- \cap \O_{+\infty}} g\,\d\mu,$
while the equality $$ \lambda \int_{\O_{-\infty} \cap
\O_{+\infty}}f_1\,\d\mu=\int_{\O_{-\infty} \cap
\O_{+\infty}}g\,\d\mu,$$ is a direct consequence of the invariance
of $\mu$ with respect to ${\mathbf{\Phi}}(\cdot,t)$. This shows that
$\lambda \|f\|_X=\|g\|_X-\|\B^+f\|_{\lp}$ for $g \geq 0$. In
general, defining
$$F_1(\x)=\int_0^{\t_-(\x)}\exp(-\lambda s)\,\left|g({\mathbf{\Phi}}(\x,-s)\right|\d s,\qquad \x \in
\O,$$ we obtain $\|\B^+f_1\|_{\lp}+\lambda \|f_1\|_X \leq
\|B^+F_1\|_{\lp}+\lambda \|F_1\|_X=\|g\|_X$ which, combined with
\eqref{B+f-2}, gives \eqref{10.1156}.
\end{proof}
\begin{nb} Notice that, in order to get the existence and uniqueness
of the solution $f$ to \eqref{BVP1}, it is not necessary for $u$ to
belong to $L^1(\Gamma_-,\d\mu_-)$. Indeed, we only have to make sure
that  $f_2 \in L^1(\O,\d\mu)$, i.e., from \eqref{f-2uu},
$\int_{\Gamma_-}|u(\y)|\left(1-e^{-\lambda
\tau_+(\y)}\right)\d\mu_-(\y) < \infty.$ Of course, to get
\eqref{10.1156}, the assumption $u \in L^1(\Gamma_-,\d\mu_-)$ is
necessary.\end{nb} Let us note that, with the notation of Theorem
\ref{Theo10.43}, we have
\begin{equation}\label{green}
\int_{\Gamma_+}\B^+f\d\mu_++\lambda
\int_{\O}f\,\d\mu=\int_{\Gamma_-}u\,\d\mu_-+\int_{\O}g\,\d\mu.
\end{equation}
Indeed, for nonnegative $u$ and $g$, \eqref{10.1156} turns out to be
precisely \eqref{green}. Then, for arbitrary $u \in \lm$ and $g \in
X$, we get \eqref{green} by splitting functions into positive and
negative parts. This leads to the following generalization of
Green's formula:
\begin{propo}[\textit{\textbf{Green's formula}}]\label{propgreen}
Let $f \in \D(\T_{\mathrm{max}})$ satisfies $\B^-f \in \lm.$ Then
$\B^+f \in \lp$ and
$$\int_{\O}\T_{\mathrm{max}}f\d\mu=\int_{\Gamma_-}\B^-f\d\mu_-
-\int_{\Gamma_+}\B^+f\,\d\mu_+$$
\end{propo}
\begin{proof} For given $f \in \D(\T_{\mathrm{max}})$, we obtain the result by setting
$u=\B^-f \in \lm$ and $g=(\lambda- \T_{\mathrm{max}})f \in X$ in Eq.
\eqref{green}.\end{proof}

\begin{nb} If $\d\mu$ is the Lebesgue measure on $\mathbb{R}^N$, the above formula leads to a better understanding of the measures $\d\mu_{\pm}$.
Indeed, comparing it to the classical Green's formula (see e.g.
\cite{bardos}), we see that the restriction of $\d\mu_{\pm}$ on the set
$\Sigma_{\pm}=\{\y \in
\partial\O\,;\,\pm \ff(\y) \cdot n(\y) > 0\}$  equals
$$\left|\ff(\y) \cdot n(\y)\right|\d\gamma(\y),$$
where $\d\gamma(\cdot)$ is the surface Lebesgue measure on
$\partial\O$.
\end{nb}

\section*{Appendix: About the class of test-functions}\setcounter{equation}{0}
\renewcommand{\theequation}{A.\arabic{equation}}

We answer in this Appendix a natural question concerning the
definition of the class of test-functions $\mathfrak{Y}$. Precisely,
we prove that two test-functions equal $\mu$--almost everywhere are
such that their derivatives (in the sense of \eqref{defi:deriva})
also coincide $\mu$-almost everywhere. To prove our claim, it
clearly suffices to prove that, given $\psi \in \mathfrak{Y}$ such
that $\psi(\x)=0$ for $\mu$-a. e. $\x \in \O$, one has
$\varphi(\x)=0$ for $\mu$-a. e. $\x \in \O$ where
$\varphi(\x)=\frac{\d}{\d s}\psi({\mathbf{\Phi}}(\x,s))\big|_{s=0}.$
Let
$$E:=\big\{\x \in \O\,;\,\psi(\x)=0 \text{ and } \varphi(\x) \neq 0\big\}.$$
It is clear that $E$ is measurable and that one has to prove that
$\mu(E)=0$. It is no loss of generality to assume that $E$
\textit{is bounded}. We observe that for any $\x \in E$, there
exists $\delta_\x
>0$ such that
\begin{equation}\label{01}\psi({\mathbf{\Phi}}(\x,t)) \neq 0, \qquad \forall\, 0 < |t| < \delta_\x.\end{equation}
Let us split $E$ as follows
$$E=\big(E \cap \O_-\big) \cup \big(E \cap \O_+ \cap \O_{-\infty}\big) \cup \big(E \cap \O_{+\infty}
\cap \O_{-\infty}\big):=E_- \cup E_+ \cup E_\infty$$ and prove that
$\mu(E_-)=\mu(E_+)= \mu(E_\infty)=0$.

\begin{enumerate}
\item First consider $E_-$. Since $\psi(\x)=0$ for $\mu$-a. e. $\x \in \O_-$ and using the fact that
any $\x \in \O_-$ can be written as  $\x={\mathbf{\Phi}}(\y,t)$ for
some $\y \in \Gamma_-$ and $0 < t < \t_+(\y)$, we observe that, for
$\mu_-$ a. e. $\y \in \Gamma_-$, $\psi({\mathbf{\Phi}}(\y,t))=0$ for
almost every (in the sense of the Lebesgue measure in $\mathbb{R}$)
$0 < t < \t_+(\y)$. For such a $\y \in \Gamma_-$, continuous
differentiability of $t \mapsto \psi({\mathbf{\Phi}}(\y,t))$ implies
$\psi({\mathbf{\Phi}}(\y,t))=0$ for {\it any } $0 < t < \t_+(\y)$.
This means, according to \eqref{01} that, for $\mu_-$-a. e. $\y \in
\Gamma_-$, ${\mathbf{\Phi}}(\y,t) \notin E$ for any $0 < t <
\t_+(\y)$. Since
$$\mu(E \cap \O_-)=\int_{\Gamma_-}\d\mu_-(\y)\int_0^{\t_+(\y)}\chi_{E}({\mathbf{\Phi}}(\y,t))\d t$$
we see that $\mu(E_-)=0.$

\item In the same way, using $\Gamma_+$ instead of
$\Gamma_-$, we show that $\mu(E \cap \O_+ \cap \O_{-\infty})=0$.

\item It remains to prove that $\mu(E_\infty)=0$. In accordance with \eqref{01},
we define for, any $n \in \mathbb{N}$,
    $$E_n:=\bigg\{\x \in E_\infty\,;\,\delta_\x \geq 1/n\bigg\}=\bigg\{\x \in E_\infty\,;\,\psi({\mathbf{\Phi}}(\x,t)) \neq 0, \:\: \forall\, 0 < |t| < 1/n\bigg\}.$$
 According to Assumption \ref{ass:h2}, it is easy to see that $\mu(E_n)=0$ for any $n \in \mathbb{N}$ since $\psi(\x)=0$ for $\mu$-a.e. $\x \in \O.$ Moreover, $E_1 \subset E_2 \subset \ldots \subset E_n \subset E_{n+1} \subset \ldots $, and
 $$\bigcap_{n=1}^\infty \bigg(E_\infty \setminus E_n\bigg)= \varnothing.$$
 Since we assumed $\mu(E) < \infty$,  we have $\mu(E_\infty \setminus E_1) < \infty$ and
 $\lim_{n \to \infty} \mu\big(E_\infty \setminus E_n\big)=0$. Writing
 $E_\infty=E_n \cup \big(E_\infty \setminus E_n\big)$, we see that $\mu(E_\infty)=0.$
\end{enumerate}

\end{document}